\documentclass[leqno]{article}
\usepackage{amssymb}
\usepackage{amsmath,  enumerate}
\usepackage{vmargin}
\usepackage[arrow,matrix]{xy}

 \makeatletter
 \renewcommand\section{\@startsection {section}{1}{\z@}%
 {-3.5ex \@plus -1ex \@minus -.2ex}%
 {2.3ex \@plus.2ex}%
 {\center \normalfont\large\bfseries}}
 \makeatother

\setmarginsrb{3cm}{3cm}{3cm}{3cm}{1mm}{6mm}{0mm}{10mm}

\numberwithin{equation}{section}

\newcommand{\8}{\infty}

\newcommand{\nz}{{\mathbb N}}
\newcommand{\nen}{n \in \nz}

\newcommand{\rz}{{\mathbb R}}

\newcommand{\qz}{{\mathbb Q}}
\newcommand{\cz}{{\mathbb C}}
\newcommand{\ten}{\otimes}
\newcommand{\tr}{{\rm tr}}
\newcommand{\un}{1\mkern -4mu{\textrm l}}
\newcommand{\ez}{{\mathbb E}\vspace{0.1cm}}

\newcommand{\Om}{\Omega}

\newcommand{\al}{\alpha}
\newcommand{\si}{\sigma}
\newcommand{\Si}{\Sigma}
\newcommand{\tet}{\theta}
\renewcommand{\d}{\delta}

\newcommand{\la}{\lambda}
\newcommand{\G}{\Gamma}

\newcommand{\eps}{\varepsilon}
\newcommand{\e}{\varepsilon}
\newcommand{\f}{\varphi}

\newcommand{\F}{{\mathcal F}}
\newcommand{\E}{{\mathcal E}}

\newcommand{\A}{{\mathcal A}}
\newcommand{\B}{{\mathcal B}}
\newcommand{\D}{{\mathcal D}}
\newcommand{\M}{{\mathcal M}}

\newcommand{\R}{{\mathcal R}}

\newcommand{\N}{{\mathcal N}}

\newcommand{\W}{{\mathcal W}}

\newcommand{\Ma}{{\mathbb M}}

\newcommand{\be}{\begin{eqnarray*}}
\newcommand{\ee}{\end{eqnarray*}}

\newcommand{\pf}{\noindent{\it Proof.~~}}
\newcommand{\qd}{\hfill$\Box$}

\newtheorem{lemma}{Lemma}[section]
\newtheorem{prop}[lemma]{Proposition}
\newtheorem{theorem}[lemma]{Theorem}
\newtheorem{cor}[lemma]{Corollary}

\newtheorem{rem}[lemma]{Remark}

\newcommand{\re}{\begin{rem}\rm}
  \newcommand{\mar}{\end{rem}}
\newtheorem{exam}[lemma]{Example}

\newcommand{\xspace}{\hbox{\kern-2.5pt}}

\newcommand{\yspace}{\hbox{\kern-0.9pt}}

\newcommand{\tnorm}[1]{\left\vert\xspace\left\vert\xspace\left\vert\mskip2mu
#1\mskip2mu \right\vert\xspace\right\vert\xspace\right\vert}

\begin{document}

\title{Noncommutative Burkholder/Rosenthal inequalities II:
applications}
\author{Marius Junge and Quanhua Xu}

\date{}

\maketitle

\begin{abstract} We show norm estimates for the sum
of independent random variables in noncommutative $L_p$-spaces for
$1<p<\infty$ following our previous work. These estimates
generalize the classical Rosenthal inequality in the commutative
case. Among applications, we derive an equivalence for the
$p$-norm of the singular values of a random  matrix with
independent entries, and characterize those symmetric subspaces
and unitary ideals which can be realized as subspaces  of a
noncommutative $L_p$ for $2<p<\infty$.
\end{abstract}


 \setcounter{section}{-1}


 \makeatletter
 \renewcommand{\@makefntext}[1]{#1}
 \makeatother \footnotetext{\noindent
 M.J.: Department of Mathematics,  University of Illinois, Urbana,
 IL 61801, USA\\
 junge@math.uiuc.edu\\
 This author is partially supported by the National Science
Foundation DMS-0301116.\\
Q.X.: Laboratoire de Math{\'e}matiques, Universit{\'e} de France-Comt{\'e},
25030 Besan\c{c}on Cedex,  France\\
qxu@univ-fcomte.fr
 }


\section{Introduction and preliminaries}


This paper is a continuation of our previous work \cite{jx-burk}
on the investigation of noncommutative martingale inequalities.
The classical theory of  martingale inequalities has a long
tradition in probability. It is well-known today that the
applications of the works of Burkholder and his collaborators
range from classical harmonic analysis to stochastic differential
equations and the geometry of Banach spaces. When proving the
estimates for the conditioned (or little) square function (cf.
\cite{burk-dis, burk-gundy}), Burkholder was aware of Rosenthal's
result \cite{ros-ineq} on sums of independent random variables.
Here we proceed differently and prove the noncommutative Rosenthal
inequality along the same line as the noncommutative Burkholder
inequality from \cite{jx-burk}. This slightly modified proof
yields a better constant. The main intention of this paper is to
illustrate the usefulness of the conditioned square function by
several examples. For many applications it is important to
consider generalized notions of independence. This will allow us
to explore applications towards random matrices and symmetric
subspaces of noncommutative $L_p$-spaces.

\smallskip

Our estimates on random matrices are motivated by the following
 noncommutative Khintchine inequality of Lust-Piquard
\cite{lust-khin}. Let $(\e_{ij})$ be an independent Rademacher
family on a probability space $(\Om, \mu)$ and let $(e_{ij})$ be
the canonical matrix units of $B(\ell_2)$. Then for any $2\le
p<\infty$ there exists a positive constant $c_p$, depending only
on $p$, such that for scalar coefficients $(a_{ij})$
 \begin{align*}
 \ez\,\big\|\sum_{ij}\e_{ij}\,a_{ij}\,e_{ij}
  \big\|_{S_p} \sim_{c_p}
 \max\Big\{\big(\sum_{i}\big(\sum_j|a_{ij}|^2 \big)^{p/2}
 \big)^{1/p}\,,\;
    \big(\sum_{j}\big(\sum_i|a_{ij}|^2
 \big)^{p/2}\big)^{1/p}\Big\}\,,
 \end{align*}
where $S_p$ denotes the usual Schatten $p$-class. Recall that for
a matrix $a=(a_{ij})$
 \[\|a\|_{S_p}= \big[\sum_n
 \la_n(|a|)^p \big]^{1/p} \,, \]
where the $\la_n(|a|)$ are the eigenvalues of $|a|$, arranged in
decreasing order and counted according to their multiplicities. In
the noncommutative setting it is natural to replace $(\e_{ij})$ by
a noncommutative independent family and the scalar coefficients
$a_{ij}$ by operator coefficients. Here we just mention, for
illustration, the following special case and refer to section
\ref{p<2} for more information. Let $(f_{ij})\subset L_p(\Om,\mu)$
be a matrix of independent mean zero random variables. Then for
$2\le p<\8$
 \begin{align*}
  &\big\|\sum_{ij}f_{ij}\ten e_{ij}
  \big\|_{L_p(\Om;S_p)}\sim_{c\,p}\\
  &~~~~ \max\Big\{\big(\sum_{ij}\|f_{ij}\|_p^p\big)^{1/p}\,,\;
  \big(\sum_{i}\big(\sum_j\|f_{ij}\|_2^2\big)^{p/2}\big)^{1/p}\,,\;
   \big(\sum_{j}\big(\sum_i\|f_{ij}\|_2^2\big)^{p/2}\big)^{1/p}
    \Big\}
   \end{align*}
and for $p<2$ (with $p'$ denoting the conjugate index of $p$)
 \begin{align*}
  &\big\|\sum_{ij}f_{ij}\ten e_{ij}
  \big\|_{L_p(\Om;S_p)}\sim_{c\,p'}\\
  &~~~~\inf\Big\{\big(\sum_{ij}\|d_{ij}\|_p^p\big)^{1/p}+
  \big(\sum_{i}\big(\sum_j\|g_{ij}\|_2^2\big)^{p/2}\big)^{1/p}+
   \big(\sum_{j}\big(\sum_i\|h_{ij}\|_2^2\big)^{p/2}\big)^{1/p}
    \Big\},
   \end{align*}
where the infimum is taken over all decompositions
$f_{ij}=d_{ij}+g_{ij}+h_{ij}$ with mean zero variables $d_{ij}$,
$g_{ij}$ and $h_{ij}$, which, for each couple $(i,j)$, are
measurable with respect to the $\si$-algebra generated by
$f_{ij}$.

The equivalence above for $p\ge 2$ is a direct consequence of our
noncommutative Rosenthal inequality in section \ref{p>2}. As
usual, the case $p<2$ is dealt with by duality. Sections \ref{p>2}
and \ref{p<2} are devoted to the Rosenthal inequalities for
$p\ge2$ and $p<2$, respectively. The random variables we consider
are general independent variables in noncommutative $L_p$-spaces
(including the type III case). In contrast with the classical case
where there exist a unique independence, one has several different
notions of independence in the noncommutative setting. Introduced
in section \ref{independence}, our definition of independence
embraces the most commonly used noncommutative notions of
independence. These include the usual tensor independence and
Voiculescu's freeness.

\smallskip

In the light of the recent concept of noncommutative maximal
functions, it would be desirable to have a perfect noncommutative
analogue of the classical Burkholder inequality by replacing the
diagonal term  $\|(d_k)\|_{\ell_p(L_p)}$ by the maximal term
$\|(d_k)\|_{L_p(\ell_\8)}$. This is indeed possible. We will make
up for it in section \ref{maximal function variant}. The same
variant is, of course, true for the noncommutative Rosenthal
inequality.

\smallskip

Symmetric subspaces of  $L_p$-spaces are motivated by
probabilistic notions of exchangeable random variables. In the
commutative situation, the memoir  of Johnson, Maurey, Schechtman
and Tzafriri \cite{jmst} contains an impressive amount of
information and many sophistical applications of probabilistic
techniques. As applications of the noncommutative
Burkholder/Rosenthal inequalities, we will extend some of their
results to the noncommutative setting in section \ref{symmetric
subspaces}. Below is an elementary example. Let $\A$ and $\M$ be
von Neumann algebras and $2\le p<\infty$. Let $(x_k)_{1\le k\le
n}\subset L_p(\M)$  and $\la>0$ such that
 \[
 \big\|\sum_{k=1}^n \e_k a_{\pi(k)} \ten  x_k\big\|_p \le\la\,
 \big\|\sum_{k=1}^n a_{k} \ten  x_k\big\|_p \]
holds for all $\e_k=\pm 1$, all  permutations $\pi$ on
$\{1,...,n\}$ and coefficients $a_k\in L_p(\A)$. Then there are
constants $\al,\beta$ and $\gamma$, depending only on $(x_k)$,
such that for all $a_k\in L_p(\A)$
 \[
 \big\|\sum_{k=1}^n a_{k} \ten  e_k\big\|_p\sim_{c_{p,\la}}
 \max\Big\{\al\big(\sum_{k=1}^n \|a_k\|_p^p\big)^{1/p}\,,\;
 \beta \big\|(\sum_{k=1}^n a_k^*a_k)^{1/2}\big\|_p\,,\;
 \gamma \big\|(\sum_{k=1}^n a_ka_k^*)^{1/2}\big\|_p\Big\}\,.
 \]
As a consequence of this statement (with $\A=\cz$), we deduce that
$\ell_p$ and $\ell_2$ are the only Banach spaces with a symmetric
basis embedding into a noncommutative $L_p$ for $2<p<\8$. On the
other hand, at the operator space level, we have four spaces
$\ell_p$, $C_p$, $R_p$ and $C_p\cap R_p$, where $C_p$ and $R_p$
are respectively the column and row subspaces of $S_p$. In the
same spirit, we also characterize the unitary ideals isomorphic to
subspaces of a noncommutative $L_p$ for $2<p<\8$ in section
\ref{unitary ideal}.

\smallskip

In the remainder of this introduction  we give some necessary
preliminaries and notation.  We use standard notation from  von
Neumann algebra theory (see e.g. \cite{kar-II, tak-I, stratila}).
For noncommutative $L_p$-spaces we follow the notation system of
\cite{jx-burk}, and refer there for more details and all
unexplained notions, especially those on martingales. As in
\cite{jx-burk}, the noncommutative $L_p$-spaces used in this paper
are those constructed by Haagerup \cite{haag-Lp}. We will work
under the standard assumptions from \cite{jx-burk}. In particular,
$\M$ is a $\si$-finite von Neumann algebra equipped with a normal
faithful state $\f$. The Haagerup noncommutative $L_p$-spaces
associated with $(\M,\f)$ are denoted by $L_p(\M)$. We denote by
$D$ the density of $\f$ in the space $L_1(\M)$ such that
 \[ \f(x)=\tr(xD), \quad x\in\M,\]
where $\tr: L_1(\M)\to\cz$ is the distinguished tracial
functional. The norm of $L_p(\M)$ is denoted by $\|\, \|_p$.
Recall that $\M D^{1/p}$ is dense in $L_p(\M)$ for any $0<p<\8$.
More generally, $D^{(1-\tet)/p}\M_a D^{\tet/p}$ is also dense in
$L_p(\M)$ for any $0\le\tet\le1$, where $\M_a$ denotes the family
of all analytic elements with respect to the modular group
$\si_t^\f$ of $\f$.

Let $\N$ be a von Neumann subalgebra of $\M$ (i.e., a w*-closed
involutive subalgebra containing the unit of $\M$). We say that
$\N$ is \emph{$\f$-invariant} if $\si_{t}^{\f}(\N)\subset\N$ for
all $t\in\rz$. According to Takesaki  \cite{tak-cond}, there
exists a unique normal faithful conditional expectation $\E:\M\to
\N$ such that $\f\circ \E=\f$. Recall that $\E$ is characterized
by
 \[\f(\E(x)y)= \f(xy)\,, \quad x\in \M,\; y\in \N.\]
Note that $\E$ commutes with the modular group $\si_t^\f$ of $\f$.
Namely, $\si_{t}^{\f}\circ \E=\E \circ \si_t^{\f}$. In these
circumstances,  $\si_{t}^{\f}\big|_\N$ is the modular group of
$\f\big|_{\N}$, and the noncommutative $L_p(\N)$ associated to
$(\N,\,\f\big|_{\N})$ can be naturally isometrically identified
with a subspace of $L_p(\M)$. With this identification, the
density of $\f\big|_{\N}$ in $L_1(\N)$ coincides with $D$. All
these allow us to not distinguish $\f$, $\si_{t}^{\f}$ and $D$ and
their respective restrictions to $\N$.

For $1\le p< \infty$, the conditional expectation $\E$ extends to
a contractive projection $\E_p$ from $L_p(\M)$ onto $L_p(\N)$
densely defined by
 \[ \E_p(xD^{1/p})= \E(x)D^{1/p}\,,\quad x\in\M.\]
$\E_p$ is also determined by
 \[\E_p(D^{(1-\theta)/p}xD^{\theta/p})=
 D^{(1-\theta)/p}\E(x) D^{\theta/p}\,, \quad
 x\in\M_a\,, \; 0\le \theta \le 1\,.\]
It is convenient to drop the index $p$. This is also justified by
using Kosaki's embedding $I:L_p(\M)\to L_1(\M)$, $I(xD^{1/p})= xD$
since then $\E_1(I(y))=I(\E_p(y))$. In this sense all maps $\E_p$
are induced by  the same map $\E_1$.

Recall that if $\N=\cz$, then $\E(x)=\f(x)1$ for every $x\in\M$;
so $\E$ can be identified with $\f$. The action of $\E$ on
$L_p(\M)$ is then given by $\E(x)=\tr(xD^{1/p'})D^{1/p}$, where
$p'$ denotes the conjugate index of $p$. Thus if additionally $\f$
is tracial, we still have $\E(x)=\f(x)1$ for $x\in L_p(\M)$.

\smallskip

We will frequently use the column, row spaces and their
conditional versions. Recall that for a finite sequence
$a=(a_k)\subset L_p(\M)$
 \[ \big\|a\big\|_{L^p(\M;\ell^c_2)}=\big\|
 \big(\sum_{k}|a_k|^2\big)^{1/2}\big\|_p\quad\mbox{and}\quad
 \big\|a\big\|_{L_p(\M;\ell^r_2)}=\big\|
 \big(\sum_{k}|a_k^*|^2\big)^{1/2}\big\|_p\,.\]
Then $L_p(\M;\ell^c_2)$ and $L_p(\M;\ell^r_2)$ are the completions
of the family of all finite sequences in $L_p(\M)$ with respect to
$\|\, \|_{L_p(\M;\ell^c_2)}$ and $\|\, \|_{L_p(\M;\ell^r_2)}$,
respectively (in the w*-topology for $p=\8$). It is convenient to
view $L_p(\M;\ell^c_2)$ and $L_p(\M;\ell^r_2)$ as the first column
and row subspaces of $L_p(B(\ell_2)\bar\otimes \M)$, respectively.

Now let $\N$ be a  $\f$-invariant von Neumann subalgebra of $\M$
with conditional expectation $\E$. Let $p\ge2$ and $a=(a_k)\subset
L_p(\M)$ be a finite sequence. Since $a_k^*a_k\in L_{p/2}(\M)$ and
$p/2\ge1$, $\E(a_k^*a_k)$ is well-defined; so we can consider
 \[\big\|a\big\|_{L_p(\M,\E;\ell^c_2)}=\big\|
 \big(\sum_{k}\E(a_k^*a_k)\big)^{1/2}\big\|_p\,.\]
According to \cite{ju-doob} (see also \cite{jx-burk}), this
defines a norm on the family of all finite sequences in $L_p(\M)$.
The corresponding completion (relative to the w*-topology for
$p=\8$) is the conditional column space $L_p(\M,\E;\ell^c_2)$.
Note that if $2\le p<\8$, then finite sequences in $\M_a D^{1/p}$
are dense in $L_p(\M,\E;\ell^c_2)$. The latter density allows us
to extend the definition to the range $1\le p<2$. Let
$a=(a_k)\subset \M D^{1/p}$ with $a_k=b_k D^{1/p}$, $b_k\in\M$.
Set
 \[\big\|a\big\|_{L_p(\M,\E;\ell^c_2)}=\big\|
 \big(\sum_{k}D^{1/p}\,\E(b_k^*b_k)D^{1/p}\big)^{1/2}\big\|_p\,.\]
We have again a norm. The resulting completion is denoted by
$L_p(\M,\E;\ell^c_2)$. The conditional row space
$L_p(\M,\E;\ell^r_2)$ is defined as the space of all $(a_k)$ such
that $(a_k^*)\in L_p(\M,\E;\ell^c_2)$, equipped with the norm
 \[\big\|(a_k)\big\|_{L_p(\M,\E;\ell^r_2)}
 =\big\|(a_k^*)\big\|_{L_p(\M,\E;\ell^c_2)}\,.\]
The space  $L_p(\M,\E;\ell^c_2)$ (resp. $L_p(\M,\E;\ell^r_2)$) can
be equally viewed as the first column (resp. row) subspace of
$L_p(B(\ell_2(\nz^2))\bar\otimes\M)$, indexed by a double index.

\begin{lemma}\label{duality of cond col Lp}
 Let $1\le p<\8$ and $p'$ be the index conjugate to $p$.
Then
 \[L_p(\M,\E;\ell^c_2)^*=L_{p'}(\M,\E;\ell^c_2)\]
holds isometrically with respect to the antilinear duality
bracket:
 \[\langle a, b\rangle=\sum\tr(b_k^*a_k)\,,\quad
 a\in L_p(\M,\E;\ell^c_2),\; b\in L_{p'}(\M,\E;\ell^c_2).\]
A similar statement holds for the conditional row spaces.
 \end{lemma}

\pf This is the column (or row) space version of \cite[Corollary
2.12]{ju-doob}. The proof there can be adapted to the present
situation by considering $\M\bar\ten B(\ell_2)$ and $\N\bar\ten
B(\ell_2)$ in place of $\M$ and $\N$, respectively. It then
remains to note that the column space $L_p(\M;\ell_2^c)$ is
complemented in $L_p(B(\ell_2)\bar\ten \M)$. See also the proof of
\cite[Theorem 2.13]{ju-doob}, where instead of one conditional
expectation, a sequence of conditional expectations is involved
(then the noncommutative Stein inequality is needed). We omit the
details.\qd

\smallskip

The preceding notations will be kept in the remainder of the
paper. Unless explicitly stated otherwise, $\M$ will denote a von
Neumann algebra equipped with a normal faithful state $\f$. If
$\N$ is a $\f$-invariant von Neumann subalgebra of $\M$, its
associated conditional expectation will be often denoted by
$\E_\N$ or simply by $\E$ if no confusion can occur.

\smallskip

The first version of this paper was written up immediately after
the submission of \cite{jx-burk} (so more than five years ago).
Since then considerable progress has been made on noncommutative
martingale inequalities. We mention only \cite{jx-const,
parcet-ran-gundy, ran-mtrans, ran-weak, ran-burk}, where, among
many other results, the optimal orders of the best constants in
most noncommutative martingale inequalities are determined.


\section{Independence}
 \label{independence}


In this section, we first introduce the central notion for our
formulation of the noncommutative Rosenthal inequality, i.e., the
independence. We then present some natural examples of
noncommutative independent variables. Our setup is the following:
$\N$ and $\A_k$ are $\f$-invariant von Neumann subalgebras of $\M$
such that $\N\subset\A_k$ for every $k$. The sequence $(\A_k)$ can
be finite.
 \begin{enumerate} [(I)]
 \item  We  say that $(\A_k)$  are (\emph{faithfully})
\emph{independent over $\N$} or \emph{with respect to $\E_\N$} if
for every $k$, $\E_\N(xy)=\E_\N(x)\E_\N(y)$ holds for all $x\in
\A_k$ and $y$ in the von Nuemann subalgebra generated by
$(\A_j)_{j\neq k}$.
 \item We  say that  $(\A_k)$ are (\emph{faithfully})
\emph{order independent} over $\N$ or \emph{with respect to
$\E_\N$} if for every $k\ge2$, $\E_{VN(\A_1, ..., \A_{k-1})}(x) =
\E_\N(x)$ holds for all $x\in \A_k$, where $VN(\A_1, ...,
\A_{k-1})$ denotes the von Neumann subalgebra generated by $\A_1,
..., \A_{k-1}$.
 \item A sequence $(x_k)\subset L_p(\M)$ is said to be \emph{faithfully}
(\emph{order}) \emph{independent with respect to $\E_\N$} if there
exist $\A_k$ such that $x_k\in L_p(\A_k)$ and $(\A_k)$ is
faithfully (order) independent with respect to $\E_\N$.
 \end{enumerate}

Note that the subalgebra $VN(\A_1, ..., \A_{k-1})$ is
$\f$-invariant too, so the conditional expectation $\E_{VN(\A_1,
..., \A_{k-1})}$ exists. Also note that the independence in (I)
can be defined for any family (without order). The adverb
\emph{faithfully} refers to the faithfulness of the state $\f$. We
will also consider the nonfaithful case in
section~\ref{nonfaithful}. If no confusion can occur, we will
often drop this adverb by saying simply independent or order
independent. If $\N=\cz$, these notions are, of course, with
respect to the state $\f$

\begin{rem}\label{independence-mart difference}
 {\rm
 Let $(\A_k)$ be order independent over $\N$. Then for every $k$
 \[\E_{VN(\A_1, ..., \A_{k-1})}(x) = \E_N(x),\quad x\in
 \A_j,\ j\ge k.\]
Indeed, we have
 \be
 \E_{VN(\A_1, ..., \A_{k-1})}(x)
 &=&\E_{VN(\A_1, ..., \A_{k-1})}\big(\E_{VN(\A_1, ...,
 \A_{j-1})}(x)\big)\\
 &=&\E_{VN(\A_1, ..., \A_{k-1})}\big(\E_\N(x)\big)
 =\E_\N(x).
 \ee
It follows that if $x_k\in L_p(\A_k)$ with $\E_\N(x_k)=0$, then
$(x_k)$ is a martingale difference sequence with respect to the
filtration $\big(VN(\A_1, ..., \A_{k})\big)_{k\ge1}$.
  }\end{rem}

\begin{lemma}\label{unc}
 Assume that $(\A_k)$ is independent over $\N$.
 \begin{enumerate}[\rm (i)]
 \item $(\A_k)$ is order independent over $\N$.
 \item If  $x_k\in L_p(\A_k)$ satisfy $\E_N(x_k)=0$, $1\le p\le\8$, then
 \[ \|\sum_{k=1}^n \e_k x_k\|_p \le 2 \|\sum_{k=1}^n
 x_k\|_p\,,\quad  \e_k=\pm 1.\]
 \end{enumerate}
\end{lemma}

\pf Let $S$ be a subset of indices and $\B_S=VN\{\A_j : j\in S\}$.
Set $\E_S=\E_{\B_S}$. Fix $k\notin S$. Now let $x\in \A_k$. We
want to prove  $\E_S(x)=\E_\N(x)$. For this  it suffices to show
 \[ \f(xy)= \f(\E_\N(x)y),\quad y\in\B_S\,.\]
This equality immediately follows from the independence of
$(\A_k)$ over $\N$ for
 \[ \f(xy)= \f(\E_\N(xy))= \f(\E_\N(x)\E_\N(y))=
 \f(\E_\N(\E_\N(x)y))=\f(\E_\N(x)y).
 \]
If we apply this to the subset $S=\{1,...,k-1\}$, we obtain (i).
To prove the second assertion consider $\eps_k=\pm 1$ and define
$S=\{k: \eps_k=1\}$. By approximation by elements of the form
$x_k=a_kD^{\frac1p}$, $a_k\in \A_k$ and $\E_N(a_k)=0$, we see that
 \[ \E_S(\sum_{k=1}^n x_k)
 =\sum_{k\in S} x_k+ \sum_{k\notin S} \E_S(x_k)
 =\sum_{k\in S} x_k .\]
Since $\E_S$ is a contraction on $L_p(\M)$,
 \[\|\sum_{k\in S} x_k\|_p=\|\E_S(\sum_{k=1}^n x_k)\|_p
 \le\|\sum_{k=1}^n x_k\|_p \,;\]
whence
 \[\|\sum_{k=1}^n \eps_k x_k\|_p
 \le \|\sum_{k\in S} x_k\|_p  +
  \|\sum_{k\in S^c} x_k\|_p
  \le 2 \|\sum_{k=1}^n x_k\|_p \,. \]
 \qd

\medskip

In the rest of this section we give some natural examples of
independent variables, which  often occur in noncommutative
probability.

\begin{exam}{\rm {\it Classical independence.} Let $(\Om, \mu)$ be a
probability space, and let $(\N, \psi)$  be a von Neumann algebra
equipped with a normal faithful state $\psi$. Let $\M=
L_\8(\Om)\bar\otimes\N$ be the von Neumann algebra tensor product
equipped with the tensor product state $\f=\mu\otimes \psi$. We
view $\N$ as a subalgebra of $\M$ in the natural way. Then the
conditional expectation $\E_\N$ is given by
 \[\E_\N(x)=\int_{\Om}xd\mu\,,\quad x\in\M\,,\]
where the integral is taken with respect to the w*-topology of
$\M$. Also recall that the noncommutative $L_p$-space $L_p(\M)$
coincides with the usual $L_p$-space $L_p(\Om;L_p(\N))$ of
$p$-integrable functions on $\Om$ with values in $L_p(\N)$. In
this case, the independence with respect to $\E_\N$ coincides with
the classical independence of vector-valued random variables. In
particular, if $(f_n)\subset L_p(\Om)$ is an independent sequence
of random variables in the usual sense, then $(f_na_n)$ is
independent with respect to $\E_\N$ for any $(a_n)\subset
L_p(\N)$.
 }\end{exam}

\begin{exam}{\rm {\it Tensor independence.} This independence is
the most transparent generalization of the classical one to the
noncommutative setting. Let $(\A_k, \f_k)$ be a sequence of von
Neumann algebras equipped with normal faithful states $\f_k$. Let
 $$(\M, \f)= \overline{\mathop{\bigotimes_{k\ge 0}}}\,(\A_k,\f_k)$$
denote the corresponding von Neumann algebra tensor product.  As
usual, we regard $\A_k$ as von Neumann subalgebras of $\M$. It is
clear that they are $\f$-invariant. The conditional expectation
$\E_{\A_k}$ is uniquely determined by
 $$\E_{\A_k}(a_0\otimes \cdots \ten a_m)
  =\big[\prod_{j\neq k} \f_j(a_j)\big]\, a_k,\quad  m\ge 0.$$
Clearly, $(\A_k)_{k\ge 1}$ is independent over  $\A_0$. If all
$\A_k$ are commutative, we go back to the classical case.
 }\end{exam}

\begin{exam}{\rm {\it Free independence.}
Our reference for this example is \cite{VDN}. Let $(\A_k)_{k\ge1}$
be a sequence of von Neumann subalgebras of $\M$, and let $\B$ be
a common von Neumann subalgebra of the $\A_k$. Assume that there
exist normal faithful conditional expectations $\E: \M\to\B$ and
$\E_k: \A_k\to\B$. The sequence $(\A_k)_{k\ge1}$ is called free
over $\B$ if
 \[\E(x_{1}\cdots x_{k})=0\]
whenever $x_{j}\in\mathring\A_{i_j}$  and
$i_1\not=i_2\not=\cdots\not=i_k$, where $\mathring\A_k=\ker\E_k$.
If $\B=\cz$, we get the freeness with respect to the state
$\f\sim\E$. There exists an equivalent way of formulating freeness
by using reduced free product.  Without loss of generality we may
assume that $\M$ is generated by the $\A_k$. Then $(\M, \E)$ can
be identified with the von Neumann algebra amalgamated  reduced
free product of  the $(\A_k, \E_k)\,$:
 $$(\M, \E)=\bar{\mathop *_{k\ge1}}\, _\B\;(\A_k,\E_k).$$
Assume in addition that $\B$ is $\si$-finite, and fix a normal
faithful state $\phi$ on $\B$. Then $\f=\phi\circ\E$  is a normal
faithful state on $\M$ and the $\A_k$ are $\f$-invariant. One
easily checks that freeness implies the independence in our sense.

Let us consider the particularly interesting case where all $\A_k$
are equal to $L_\8(-2,2)$, equipped with the Wigner measure
 $$d\mu(t)=\frac{1}{2\pi}\,\sqrt{4-t^2}\,dt.$$
Then the reduced free product (without amalgamation)
 $$(\M, \f)=\bar{\mathop *_{k\ge1}}\;\A_k$$
is a II$_1$ factor with $\f$ a normal faithful tracial state. Let
$x_k\in\A_k$ be given by $x_k(t)=t$. Then the sequence $(x_k)$ is
free. This is a semicircular system in Voiculescu's sense. It is
the free analogue of a standard Gaussian system.

\smallskip

Semicircular systems admit a more convenient realization via Fock
spaces. Let us describe this briefly. Let $H$ be a complex Hilbert
space. The associated  free (or full) Fock space  is defined by
 $$\F(H) = \bigoplus_{n \geq 0}H^{\otimes n},$$
where $H^{\ten 0} = {\cz}\un$ ($\un$ being a unit vector, called
vacuum), and $H^{\ten n}$ is the n-th Hilbertian tensor power of
$H$ for $n\ge1$. The (left) creator associated with a vector
$\xi\in H$ is the operator on $\F(H)$ uniquely determined by
 $$c(\xi)\,\xi_1\ten\cdots \ten\xi_n
 =\xi\ten\xi_1\ten\cdots\ten \xi_n$$
for any $\xi_1, ..., \xi_n\in H$. Here $\xi_1\ten\cdots\ten \xi_n$
is understood as the vacuum $\un$ if $n=0$. Its adjoint is given
by
 $$c(\xi)^*\,\xi_1\ten\cdots\ten \xi_n=\langle\xi_1,\;\xi\rangle\,
 \xi_2\ten\cdots\ten \xi_n$$
(with $c(\xi)^*\un=0$). This is the annihilator associated with
$\xi$ and is denoted by $a(\xi)$.  We have the following free
commutation relation:
 $$ a(\eta)c(\xi)=\langle\xi,\;\eta\rangle 1,
 \quad\xi,\;\eta\in H.$$
Now assume that $H$ is the complexification of a real Hilbert
space $H_{\rz}$.  For a real $\xi\in H_{\rz}$  define
 $$g(\xi)=c(\xi)+a(\xi).$$
Let  $\G(H)$ be the von Neumann subalgebra of $B(\F(H))$ generated
by all $g(\xi)$ with real $\xi\in H_{\rz}$:
 $$\G(H)=\big\{g(\xi): \xi\in H_{\rz}\big\}''\,.$$
This is  the  free von Neumann algebra associated with $H$ (or
more precisely, with $H_{\rz}$). The vector state $\f$ defined by
the vacuum, $x\mapsto \langle x\un,\; \un\rangle$ is faithful and
tracial on $\G(H)$. If $(\xi_k)$ is an orthonormal system of $H$
consisting of real vectors, then $(g(\xi_k))$ is a semicircular
system.

\smallskip

The preceding Fock space construction can be deformed to get type
III algebras. For this let $H$ be separable and fix an orthonormal
basis $(e_{\pm k})_{k\ge1}$ of $H$ consisting of real vectors. Let
$\la=(\la_k)$ be a sequence of positive numbers. Set
 \begin{equation}\label{shlya}
 g_k=c(e_k)+\sqrt{\la_k}\; a(e_{-k})\,,\quad k\ge1\,.
 \end{equation}
Let $\G_\la$ be the von Neumann algebra on $\F(H)$ generated by
$(g_k)$, and let $\f_\la$ be the vector state on $\G_\la$
determined by the vacuum. Then $(g_k)$ is free in $(\G_\la,
\f_\la)$. This is a generalized circular system in Shlyakhtenko's
sense \cite{shlya-quasifree}. If all $\la_k$ are equal to $1$,
$\G_\la$ becomes the previous free von Neumann algebra $\G(H)$
associated with $H$. Otherwise, $\G_\la$ is a type III factor and
the state $\f_\la$ is called a free quasi-free state.
 }\end{exam}

\begin{exam}{\rm {\it $q$-independence.} The Fock space
construction in the previous example can be modified to embrace
the so-called $q$-independence, $-1\le q\le 1$, introduced by
Bo\.zejko and Speicher \cite{bos-example, bos-cox, boks}. Again,
let $H$ be the complexification of a real Hilbert space $H_{\rz}$.
The associated $q$-Fock space ${\F}_{q}(H)$ is defined by
 $${\F}_q(H) = \bigoplus_{n \ge 0}H^{\ten n},$$
where  $H^{\ten n}$ is now equipped with the $q$-scalar product
for every $n\ge2$. Recall that ${\F}_0(H)$ is the free Fock space
discussed in the previous example, while ${\F}_1(H)$ and
${\F}_{-1}(H)$ are the classical symmetric and  antisymmetric Fock
spaces, respectively.

Given $\xi\in H$  we define the corresponding creator $c_q(\xi)$
and annihilator $a_q(\xi)$ similarly as in the free case. These
are linear operators on $\F_q(H)$ determined by the following
conditions
 $$c_q(\xi)\,\xi_{1}\ten\cdots\ten \xi_{n}
 =\xi\ten \xi_{1}\ten\cdots\ten\xi_{n}$$
and
 $$a(\xi)\,\xi_{1}\ten\cdots\ten\xi_{n}
 =\sum_{k=1}^nq^{k-1}\langle \xi_k,\; \xi\rangle\,
 \xi_1\ten\cdots\ten{\mathop \xi^\vee}_k\ten\cdots\ten \xi_{n},$$
where ${\displaystyle\mathop \xi^\vee}_k$ means that $\xi_k$ is
removed from the tensor product. $c_q(\xi)$ and $a_q(\xi)$ are
bounded operators if $q<1$ and closable densely defined operators
if $q=1$. In the latter case, $c_q(\xi)$ and $a_q(\xi)$ also
denote their closures. Again, we have $c_q(\xi)^*=a_q(\xi)$. The
creators and annihilators satisfy the following  $q$-commutation
relations :
 \[a_q(\xi)c_q(\eta)-q\,c_q(\eta)a_q(\xi)
 =\langle \eta,\; \xi\rangle 1,\quad \xi,\;\eta\in H.\]
In the cases of $q=\pm1$ these are respectively the canonical
commutation relations (CCR) and the  canonical anticommutation
relations (CAR).

Given a real vector $\xi\in H_\rz$ define
 $$g_q(\xi)=c_q(\xi)+a_q(\xi).$$
$g_q(\xi)$ is called a $q$-Gaussian variable. The $q$-von Neumann
algebra $\G_{q}(H)$ associated with $H$ is the von Neumann algebra
on $\F_q(H)$ generated by the $g_q(\xi)$ with real $\xi$. As in
the free case, the vacuum expectation $x\mapsto \langle x\un,
\un\rangle$ is a normal faithful tracial state on $\G_q(H)$,
denoted by $\tau_q$. In particular, $\G_{0}(H)$ is the free  von
Neumann algebra considered previously. On the other hand, if $\xi$
and $\eta$ are orthogonal, then $g_1(\xi)$ and $g_1(\eta)$
commute, while $g_{-1}(\xi)$ and $g_{-1}(\eta)$ anticommute.
Therefore, $\G_1(H)$ is commutative, while $\G_{-1}(H)$ is a
Clifford algebra.

Let $K\subset H$ be a closed subspace, which is the
complexification of $K_\rz\subset H_\rz$. Then $\G_q(K)$ is a
subalgebra of $\G_q(H)$. The associated conditional expectation is
given by the second quantization of the orthogonal projection from
$H_\rz$ onto $K_\rz$.  Now let $(H_k)$ be a sequence of subspaces
of $H$ which are complexifications of pairwise orthogonal
subspaces of $H_\rz$. Each $\G_q(H_k)$ is identified with the von
Neumann subalgebra of $\G_q(H)$ generated by $g_q(\xi)$ with real
$\xi\in H_k$. Then the $\G_q(H_k)$ are independent with respect to
$\tau_q$. Consequently, if $(\xi_k)_{k}$ is an orthonormal
sequence of real vectors of $H$, $(g_q(\xi_k))_k$ is independent.
This sequence $(g_q(\xi_k))_k$ is called a $q$-semicircular
system.

\smallskip

Shlyakhtenko's generalized circular systems admit $q$-counterparts
too. We refer to \cite{hiai} for more details. Here we briefly
discuss only the case $q=-1$, which is a reformulation of the
classical construction of the Araki-Woods factors. These latter
factors are built using Pauli matrices as follows. We consider the
generators of the CAR algebra
 \begin{equation}\label{aw}
 a_k =1 \ten \cdots \ten 1  \ten
 \underbrace{e_{12}}_{\mbox{\scriptsize $k$-th position}}\ten
 1\ten \cdots \ten 1
 \end{equation}
in the algebraic tensor product $\otimes_{k\ge 1}\,
 \Ma_2$, where, as usual, $e_{ij}$ denote
the matrix units of $\Ma_2=B(\ell_2^2)$. Fix a sequence
$(\mu_k)\subset (0,1)$, and consider the states
$\f_k=(1-\mu_k)e_{11} +\mu_ke_{22}$ on $\Ma_2$. Then the tensor
product state $\f=\ten_{k\ge1} \f_{k}$ is a quasi-free state
satisfying
 \[ \f(a_{i_1}^*\cdots a_{i_r}^*\,a_{j_1}\cdots a_{j_s})=
 \delta_{rs} \prod_{l=1}^s \delta_{i_l,j_l}\,\mu_{i_l}\]
for all increasing sequences $i_1<...<i_r$ and $j_1<...<j_s$. We
denote by $\W$ the von Neumann algebra generated by the $a_k$'s in
the GNS construction with respect to $\f$. Then $\W$ is a
hyperfinite type III factor and $(a_k)$ are independent with
respect to $\f$.

 }\end{exam}

\begin{exam}{\rm {\it Group algebras.} Consider a discrete group $G$. Let
$VN(G)\subset B(\ell_2(G))$ be the associated von Neumann algebra
generated by the left regular representation $\la: G\to
B(\ell_2(G))$. More precisely, $\la$ is defined by
 \[\big(\la(g)f\big)(h)=f(g^{-1}h),\quad f\in\ell_2(G),\;
 h, g\in G\,,\]
and $VN(G)$ is generated by $\{\la(g)\ :\ g\in G\}$. Recall that
$VN(G)$ is also the w*-closure in $B(\ell_2(G))$ of the algebra of
all finite sums $\sum \al(g)\la(g)$ with $\al(g)\in\cz$. Let
$\tau_{G}$ be the vector state on $VN(G)$ determined by $\d_e$,
where $e$ is the identity of $G$ and $(\d_g)_{g\in G}$ is the
canonical basis of $\ell_2(G)$.  $\tau_{G}$ is a normal faithful
tracial state on $VN(G)$. If $H$ is a subgroup of $G$, then
$VN(H)$ is identified with the von Neumann subalgebra of $VN(G)$
generated by $\{\la(h) : h\in H\}$. The corresponding conditional
expectation $\E_{VN(H)}$ is determined by
 \[\E_{VN(H)}\big[\sum_{g\in G} \al(g)\la(g)\big]
 =\sum_{g\in H}\al(g)\la(g),\quad\al(g)\in\cz\,.\]
Now let $(G_n)$ be an increasing sequence of subgroups of $G$ and
$g_n\in G_n\setminus G_{n-1}$. Then it is easy to see that
$(\la(g_n))_n$ is order independent (but not independent in
general) with respect to $\tau_G$. In particular, a sequence of
free generators on a free group is order independent. Moreover, it
is clearly independent.
 }\end{exam}


\section{Noncommutative Rosenthal inequality: $p\ge2$}
 \label{p>2}


In this section we prove the noncommutative Rosenthal inequality
in the case $p\ge2$. In this section $\M$ will denote a von
Neumann algebra with a normal faithful state $\f$, and
$\N\subset\M$ a $\f$-invariant von Neumann subalgebra with
conditional expectation $\E=\E_\N$.  Following \cite{jx-burk}, we
will also need the diagonal space $\ell_p(L_p(\M))$ whose norm
will be denoted by $\|\,\|_{\ell_p(L_p)}$. In the remainder of the
paper, $c$ will denote an absolute positive constant which may
change from line to line, and $c_p$ a positive constant depending
only on $p$. The notation $A\sim_{c}B$ will mean that $A\le c\,B$
and $B\le c\,A$.

\begin{theorem}\label{ros}
 Let $2\le p<\infty$ and $(x_k)\in L_p(\M)$ be a finite sequence
such that $\E(x_k)=0$.
 \begin{enumerate}[\rm(i)]
 \item  If $(x_k)$ is  independent with respect to $\E$,
 then
 \begin{align*}
  \frac{c}{p}\, \big\|\sum_{k} x_k \big\|_p
  &\le \max\big\{\|(x_k)\|_{\ell_p(L_p)},
  \|(x_k)\|_{L_p(\M,\E;\ell_2^c)}, \|(x_k)\|_{L_p(\M,\E;\ell_2^r)}\big\}
  \le 2 \big\|\sum_{k} x_k \big\|_p \,.
 \end{align*}
 \item  If $(x_k)$ is order independent with respect to $\E$,
 then
   \begin{align*}
  \frac{c}{p^2}\, \big\|\sum_{k} x_k \big\|_p
  &\le \max\big\{\|(x_k)\|_{\ell_p(L_p)},
  \|(x_k)\|_{L_p(\M,\E;\ell_2^c)}, \|(x_k)\|_{L_p(\M,\E;\ell_2^r)}\big\}
  \le 2 \big\|\sum_{k} x_k \big\|_p \,.
 \end{align*}
 \end{enumerate}
 \end{theorem}

\pf (i) Let $(\A_k)$ be a sequence of $\f$-invariant von Neumann
subalgebras of $\M$ which are independent over $\N$ and such that
$x_k\in L_p(\A_k)$. Then by Lemma \ref{unc} (ii) and the fact that
$L_p(\M)$ is of cotype $p$ with constant $1$, we obtain
 \[ \|(x_k)\|_{\ell_p(L_p)}\le 2\big\|\sum_{k} x_k \big\|_p\,.\]
On the other hand, by independence,
 \[\E(x^*_kx_j)=0,\quad k\neq j.\]
Thus, for $x=\sum x_k$, we have
 \begin{align*}
  \big\|\sum_{k} \E(x_k^*x_k)\big\|_{p/2}
  &=\big\|\E(x^*x)\big\|_{p/2} \le
  \big\|x^*x\big\|_{p/2}= \| x\|_p^2,.
 \end{align*}
Therefore the lower estimate for the norm of the sum is proved.

The main part is the proof of the upper estimate. First, let us
observe that this  upper estimate is  also true for $1\le p\le2$
since
 \begin{equation}\label{[1,2]}
 \big\|\sum_{k} x_k \big\|_p\le
 2 \|(x_k)\|_{\ell_p(L_p)} \,.
 \end{equation}
Indeed, this inequality follows immediately from the
unconditionality of $(x_k)$ given by Lemma \ref{unc} (ii) and the
type $p$ property of $L_p(\M)$.  To treat the case $p\ge2$ we will
use a standard iteration procedure. The key step is to show that
if the upper estimate is true for some $p\ge1$, then so is it for
$2p$. This will enable us to iterate, by using \eqref{[1,2]} as a
starting point. Thus we assume that for some $p$ there exists a
positive constant $c_p$ such that
 \[\big\|\sum_{k} x_k \big\|_p
  \le c_p \max\big\{\|(x_k)\|_{\ell_p(L_p)},\
  \|(x_k)\|_{L_p(\M,\E;\ell_2^c)},\
  \|(x_k)\|_{L_p(\M,\E;\ell_2^r)}\big\}\]
for all $x_k\in L_p(\A_k)$ with  $\E(x_k)=0$. Our aim is to prove
the same estimate for $2p$.  Let $x_k\in L_{2p}(\A_k)$ and
$\E(x_k)=0$. First, we apply the noncommutative Khintchine
inequality (cf. \cite{LPP} and also \cite{pis-ast} with the right
order of the best constant) and deduce from Lemma \ref{unc} that
 \begin{equation}\label{square}
 \big\|\sum_{k} x_k \big\|_{2p}\le
 2\,\ez \big\|\sum_{k}\eps_k x_k\big\|_{2p}\le c\sqrt{p}\,
 \max\big\{\big\|\sum_{k} x_k^*x_k\big\|_{p}^{1/2},\;
 \big\|\sum_{k} x_kx_k^*\big\|_{p}^{1/2}\big\},
 \end{equation}
where $(\eps_k)$ is a Rademacher sequence and $\ez$ denotes the
corresponding expectation. Let us consider the first square
function on the right hand side. We define the mean zero elements
$y_k=x_k^*x_k-\E(x_k^*x_k)$. By assumption, we have
 \be
 \big\|\sum_{k} x_k^*x_k\big\|_{p}
 &\le&\big\|\sum_{k}\E(x_k^*x_k)\big\|_{p}
 + \big\|\sum_{k} y_k\big\|_{p}\\
 &\le& \big\|\sum_{k}\E(x_k^*x_k)\big\|_{p}  +
 c_p\max\big\{\big\|(y_k)\big\|_{\ell_p(L_p)},\;
  \big\|(y_k)\big\|_{L_p(\M,\E;\ell_2^c)}\big\}
 \ee
Moreover, if $1\le p\le 2$,  we can disregard the second term in
the maximum by virtue of (\ref{[1,2]}). Since $\E$ is a
contraction on $L_p(\M)$, we have
 \[ \big\|(y_k)\big\|_{\ell_p(L_p)}
 = \big(\sum_{k}\big\|x_k^*x_k-\E(x_k^*x_k)\big\|_p^p \big)^{1/p}
 \le 2  \big(\sum_{k} \|x_k\|_{2p}^{2p} \big)^{1/p} \,.\]
Hence, for $1\le p\le 2$, we find
 \begin{align*}
  \big\|\sum_{k} x_k \big\|_{2p}\le
  c\sqrt{5p}\,  \max\big\{\|(x_k)\|_{\ell_p(L_p(\M))},\
  \|(x_k)\|_{L_p(\M,\E;\ell_2^c)},\ \|(x_k)\|_{L_p(\M,\E;\ell_2^r)}\big\}.
 \end{align*}
Now assume $2<p<\infty$. We first note that
 \be
 \E(y_k^2) &=&
 \E\big[\big(x_k^*x_k-\E(x_k^*x_k)\big)^*
 \big(x_k^*x_k-\E(x_k^*x_k)\big)\big]\\
 &=&
 \E(x_k^*x_kx_k^*x_k)-\E(x_k^*x_k)\E(x_k^*x_k)
 \le  \E(|x_k|^4).
 \ee
Using \cite[Lemma 5.2]{jx-burk}, we obtain
 \begin{align*}
 \big\|\sum_k\E(|x_k|^4)\big\|_{p/2}
 \le \big\|\sum_{k} \E(|x_k|^2)\big\|_p^{(p-2)/(p-1)}
 \big(\sum_{k} \|x_k\|_{2p}^{2p} \big)^{1/(p-1)} \,.
 \end{align*}
By homogeneity, this implies
 \[ \big\|\sum_k\E(|x_k|^4)\big\|_{p/2}^{1/2}\le
 \max\big\{\big\|(x_k)\big\|_{\ell_{2p}(L_{2p})}^2\,,\;
  \big\|(x_k)\big\|_{L_{2p}(\M,\E;\ell_2^c)}^2\big\}.
  \]
Therefore we have proved that
 \[ \big\|(y_k)\big\|_{L_{p}(\M,\E;\ell_2^c)}\le
 \max\big\{\big\|(x_k)\big\|_{\ell_{2p}(L_{2p})}^2\,,\;
  \big\|(x_k)\big\|_{L_{2p}(\M,\E;\ell_2^c)}^2\big\}.
  \]
Applying the same arguments to $x_kx_k^*$ and putting together all
inequalities so far obtained,  we find
 \begin{align*}
  \big\|\sum_{k} x_k \big\|_{2p}\le
 c(p(1+2c_p))^{1/2}
  \max\big\{\|(x_k)\|_{\ell_{2p}(L_{2p})},\
  \|(x_k)\|_{L_{2p}(\M,\E;\ell_2^c)},\
  \|(x_k)\|_{L_{2p}(\M,\E;\ell_2^r)}\big\}.
  \end{align*}
It thus follows that
 $$c_{2p}\le c(p(1+2c_p))^{1/2}$$
for $p>2$. We then deduce that $c_{2p}\le c'2p$ for some absolute
constant $c'$. Therefore, the induction argument works and we
obtain assertion (i).

(ii) The proof of this part is almost the same as the previous
one. The only difference is that Lemma \ref{unc} is no longer at
our disposal. In consequence, we have to replace \eqref{square} by
the noncommutative Burkholder-Gundy inequality from \cite{px-BG,
jx-burk} (see also \cite{jx-const} for the right order of the best
constants):
 \[
 \big\|\sum_{k} x_k \big\|_{2p}\le cp
 \max\big\{\big\|\sum_{k} x_k^*x_k\big\|_{p}^{1/2},\;
 \big\|\sum_{k} x_kx_k^*\big\|_{p}^{1/2}\big\}.
  \]
This is true for $(x_k)$ is a martingale difference sequence.
Indeed, since the von Neumann subalgebra generated by the $\A_k$
is $\f$-invariant, we may assume that this subalgebra is $\M$
itself. Then letting $\M_k=VN(\A_1, ..., \A_k)$, we see that
$(\M_k)$ is an increasing filtration of subalgebras in the sense
of \cite{jx-burk}, which yields a noncommutative martingale
structure in $\M$. By Remark \ref{independence-mart difference},
$(x_k)$ is a martingale difference sequence with respect to
$(\M_k)$. The rest of the proof is then the same as that of (i).
\qd

\begin{rem}{\rm
 In the commutative case the best constant in
the Rosenthal inequality is of order $p/(1+\log p)$ as $p\to\8$
(cf. \cite{jsz}). In view of this result, the constant of order
$p$ in the first inequality in Theorem \ref{ros} seems reasonable.
At the time of this writing we do not know whether this order is
optimal.
 }\end{rem}

Theorem \ref{ros} deals with independent mean zero variables.  For
general independent variables, we have the following easy
consequence. From now on we will confine our attention only to
independence. All subsequent results have counterparts for order
independence.

\begin{cor}
 Let $p$ and $\M$ be as in Theorem \ref{ros}. Let $(x_k)\subset
 L_p(\M)$ be an independent sequence with respect to $\E$. Then
  \[\|\sum_{k}x_k \|_{p}\le cp
  \max\big\{\|\sum_{k}\E(x_k)\|_{p}\,,\;\|(x_k)\|_{\ell_p(L_p)},\
  \|(x_k)\|_{L_p(\M,\E;\ell_2^c)},\ \|(x_k)\|_{L_p(\M,\E;\ell_2^r)}\big\}.
  \]
 If additionally all $x_k$ are positive, the inverse inequality
 holds without constant.
\end{cor}

\pf Let $y_k=x_k-\E(x_k)$. Then
 \[\|\sum_{k}x_k \|_{p}\le \|\sum_{k}\E(x_k)\|_{p}+
 \|\sum_{k}y_k\|_{p}\,.\]
Now applying Theorem \ref{ros} to the centered sequence $(y_k)$,
we get an equivalence for the second term on the right. Using
triangle inequality and $\|\E(x_k)\|_p\le\|x_k\|_p$, we have
 \[\|(y_k)\big\|_{\ell_p(L_p)}\le
 2\|(x_k)\big\|_{\ell_p(L_p)}\,.\]
For the terms on the conditional square functions, we note that
 \[\E(|y_k|^2)=\E(|x_k|^2) - |\E(x_k)|^2\le\E(|x_k|^2).\]
Then we deduce the desired inequality. To prove the additional
part, by the contractivity of $\E$ on $L_p(\M)$
 \[\|\sum_{k}x_k \|_{p}\ge \|\sum_{k}\E(x_k)\|_{p}\,.\]
On the other hand, by  Jensen's inequality
 \begin{equation}\label{lower khintchine}
 \big\| \sum_k |x_k|^2 \big\|_{p/2}=
 \big\|\ez\big(|\sum_k\eps_kx_k|^2\big)\big\|_{p/2}\le
 \ez\big\|\sum_k\eps_kx_k\big\|_{p}^2\,.
 \end{equation}
Note that since $x_k\ge0$, $-\sum x_k\le\sum\eps_kx_k\le\sum x_k$
for any $\eps_k=\pm1$; so $\|\sum\eps_kx_k\big\|_{p}\le \|\sum
x_k\big\|_{p}$. Therefore,
 \[\|(x_k)\|_{L_p(\M,\E;\ell_2^c)}\le
 \big\|\big(\sum_k |x_k|^2\big)^{1/2}\big\|_{p}
 \le\|\sum_k x_k\|_{p}\,.\]
For the diagonal term, it suffices to note the inequality
 \begin{equation}\label{p-2}
 \|(x_k)\|_{\ell_p(L_p)}
 \le\big\|\big(\sum_k |x_k|^2\big)^{1/2}\big\|_{p}\,,
 \end{equation}
which is obtained by interpolating the two cases $p=2$ and $p=\8$.
Thus the proof of the corollary is complete. \qd

\smallskip

In the case $\N=\cz$, our  Rosenthal inequality takes a simpler
form. Let us formulate this explicitly as follows.

\begin{cor}\label{scalar ros}
 Let $2\le p<\8$, and let $(x_k)\subset L_p(\M)$ be a  sequence
independent with respect to $\f$ such that $\tr(x_kD^{1/p'})=0$.
Then
 \begin{align*}
 \big\|\sum_{k}x_k \big\|_{p}\,\sim_{cp}
 \max\big\{\big(\sum_{k}\|x_{k}\|_p^p\big)^{1/p}\,, \;
  \big(\sum_{k}\tr[(x_{k}^*x_k +x_kx_k^*)D^{1-2/p}]\big)^{1/2}
    \big\}.
 \end{align*}
In particular, if $\f$ is tracial,
 \begin{align*}
 \big\|\sum_{k}x_k \big\|_{p}\,\sim_{cp}
 \max\big\{\big(\sum_{k}\|x_{k}\|_p^p\big)^{1/p}\,, \;
  \big(\sum_{k}\|x_{k}\|_2^2\big)^{1/2}
    \big\}.
 \end{align*}

\end{cor}

\pf It suffices to observe that for any $q\ge 1$ the conditional
expectation $\E_\cz$ on $L_q(\M)$ is given by
$\E_\cz(x)=\tr(xD^{1/q'}) D^{1/q}$. \qd

\smallskip

In the same spirit, we have the following Khintchine type
inequality.

\begin{cor}\label{khin-p2}
 Keep the assumptions of Corollary \ref{scalar ros} and assume in addition that
 \[0<\kappa_1=\inf_k \tr[(x_{k}^*x_k +x_kx_k^*)D^{1-2/p}]\quad
 \mbox{and}\quad \sup_k\|x_k\|_p=\kappa_2<\8\,.\]
Let $\A$ be another von Neumann algebra, and let $(a_k)\subset
L_p(\A)$. Then
 \[
 \big\|\sum_{k}a_k\ten x_k \big\|_{L_p(\A\bar\ten\M)}
 \,\sim_{c_{p, \kappa_1, \kappa_2}}
  \big\|\big(\sum_{k}a_k^*a_k + a_ka_k^*\big)^{1/2}\big\|_p\,.
 \]
\end{cor}

\pf We may assume that $\A$ is $\si$-finite, so equipped with a
normal faithful state $\psi$. Then the tensor product
$\A\bar\ten\M$ is equipped with $\psi\ten\f$. Identifying $\A$
with a subalgebra of $\A\bar\ten\M$ by $a\leftrightarrow a\ten1$,
we see that the associated conditional expectation satisfies
$\E_\A(a\ten x)=\tr(xD^{1/p'})a$ for $a\in L_p(\A)$ and $x\in
L_p(\M)$. The independence of $(x_k)$ with respect to $\f$ implies
that of $(a_k\ten x_k)$ with respect to $\E_\A$. Therefore, by
Theorem \ref{ros}, we obtain an equivalence of
$\big\|\sum_{k}a_k\ten x_k \big\|_p$ with the maximum of three
terms. Let us first consider the two terms on the conditional
square functions:
  \be
  \big\|(a_k\ten x_k)\big\|_{L_p(\A\bar\ten\M, \E_\A;\ell_2^c)}
  &=&\big\|\big(\sum_{k}a_k^*a_k\ten\tr(x_k^*x_kD^{1-\frac2p})
  D^{\frac2p}\big)^{1/2}\big\|_p \\
  &=&\big\|\big(\sum_{k}\tr(x_k^*x_kD^{1-\frac2p})\,a_k^*a_k
  \big)^{1/2}\big\|_p\\
  &\ge & \sqrt{\kappa_1}\,\big\|\big(\sum_{k}a_k^*a_k
  \big)^{1/2}\big\|_p \,.
  \ee
On the other hand, by the H\"older inequality,
 \[\tr(x^*x_kD^{1-2/p})\le \|x_k\|_p^2\le \kappa_2^2\,.\]
Thus it follows that
 \[\big\|(a_k\ten x_k)\big\|_{L_p(\A\bar\ten\M, \E_\A;\ell_2^c)}\sim
 \big\|\big(\sum_{k}a_k^*a_k\big)^{1/2}\big\|_p \,.\]
Passing to adjoints, we get the same estimate for the other
conditional square function. Similarly, we have
 \[\big\|(a_k\ten x_k)\big\|_{\ell_p(L_p)}\sim
 \big\|(a_k)\big\|_{\ell_p(L_p)}\,.\]
However, by (\ref{p-2})
 \[\big\|(a_k)\big\|_{\ell_p(L_p)}
 \le\big\|\big(\sum_{k}a_k^*a_k\big)^{1/2}\big\|_p \,.\]
Therefore, the assertion follows. \qd

\smallskip

We end this section by a remark on general von Neumann algebras.

\begin{rem}\label{general VN}{\rm
 As stated, our noncommutative Rosenthal inequality holds
for $\si$-finite von Neumann algebras. It can  be easily extended
to an arbitrary von Neumann algebra $\M$ provided $\N$ and
$(\A_k)$ are von Neumann subalgebras of $\M$ such that there exist
normal faithful conditional expectations $\E_\N:\M\to \N$ and
$\E_{\A_k}: \M\to\A_k$ satisfying the commutation relation
$\E_{\A_k}\E_{\N}=\E_{\N}\E_{\A_k}=\E_{\N}$. Indeed, let $\psi$ be
a strictly normal semifinite faithful weight on $\N$, i.e., a
weight of the form $\psi=\sum_{i\in I} \phi_i$, where the $\phi_i$
are normal states on $\N$ with mutually orthogonal supports. Let
$e_i$ be the support of $\phi_i$.  For a finite subset $J\subset
I$, set $e_J=\sum_{i\in J} e_i$. Then $(e_J)$ is an increasing
family of projections such that $\lim_J e_J=1$ strongly. Now we
may consider the normal faithful state
 $$\f_J=\frac1{|J|}\,\sum_{i\in J} \phi_i\circ\E_\N\quad
 \mbox{on}\quad e_J\M e_J\,.$$
If $(x_k)\subset L_p(\M)$ is an independent sequence with respect
to $\E_\N$ and $(\A_k)$ is the associated independent sequence of
subalgebras, we see that the assumptions of Theorem \ref{ros} are
satisfied for $\A_{k,J}=e_J\A_ke_J$. Moreover, for $x\in L_p(\M)$
with $p<\8$ we have
 $$x=\lim_J e_Jx=\lim_J xe_J=\lim_J e_Jxe_J\quad
 \mbox{in}\quad L_p(\M).$$
Thus  by density, Theorem \ref{ros} holds in $L_p(\M)$. This
remark applies to all results proved in this paper. We will not
repeat it and consider only the $\si$-finite case for simplicity.
 }\end{rem}


\section{Noncommutative Rosenthal  inequality: $p<2$}
 \label{p<2}


We now investigate  the noncommutative Rosenthal inequality for
$1<p\le 2$, which is the dual version of Theorem \ref{ros}. As for
the Burkholder inequality in \cite{jx-burk}, this dual version did
not exist explicitly in literature even in the commutative
(=classical) case. In this section we will assume as before that
$\N$ and $(\A_k)$ are $\f$-invariant von Neumann subalgebras of
$\M$ such that $(\A_k)$ is independent with respect to the
conditional expectation $\E=\E_\N$.

\smallskip

We start by considering the subspace $\R_p^c$ of
$L_p(\M,\E;\ell_2^c)$ consisting of all sequences $(x_k)$ such
that $x_k\in L_p(\A_k)$ with $\E(x_k)=0$, $1\le p<\8$.
Alternately, $\R_p^c$ can be defined as the closure in
$L_p(\M,\E;\ell_2^c)$ of all sequences $(a_kD^{1/p})$ such that
$a_k\in\A_k$ with $\E(a_k)=0$. Similarly, we define the
corresponding subspaces of $L_p(\M,\E;\ell_2^r)$ and
$\ell_p(L_p(\M))$, which are denoted respectively by $\R_p^r$ and
$\R_p^d$ .

\begin{lemma}\label{r-compl}
Let $1\le p < \infty$. Then $\R_p^c$ is 2-complemented in
$L_p(\M,\E;\ell_2^c)$. The similar statements hold for the row and
diagonal subspaces $\R_p^r$ and $\R_p^d$.
\end{lemma}

\pf Let us consider a finite sequence $(a_kD^{1/p})$ with
$a_k\in\M$. By the Cauchy-Schwarz inequality
 \[\E\big(\E_{\A_k}(a_k)^*\E_{\A_k}(a_k)\big)
 \le\E\big(\E_{\A_k}(a_k^*a_k)\big)=\E(a_k^*a_k).\]
It follows that
 \[\big\|(\E_{\A_k}(a_k)D^{1/p})\big\|_{L_p(\M,\E;\ell_2^c)}
 \le\big\|(a_kD^{1/p})\big\|_{L_p(\M,\E;\ell_2^c)}\,.\]
This shows that the map $F((x_k))=(\E_{\A_k}(x_k))$ defines a
contraction on $L_p(\M,\E;\ell_2^c)$. The same argument shows that
$E((x_k))=(\E(x_k))$ is also a contraction. Then $({\rm id}-E)F$
is the desired projection from $L_p(\M,\E;\ell_2^c)$ onto
 $\R_p^c$. This same projection is also bounded from
$L_p(\M,\E;\ell_2^r)$ onto $\R_p^r$ and from $\ell_p(L_p(\M))$
onto $\R_p^d$.\qd

\begin{theorem}\label{dual}
Let $1<p\le 2$. Let  $x_k\in L_p(\A_k)$ such that $\E(x_k)=0$.
Then
 \[\frac{1}{2}\, \|\sum_k x_k \|_p \le
 \inf_{x_k=x_k^d+x_k^c+x_k^r}\,
 \big\{\|(x_k^d)\|_{\R_p^d} + \|(x_k^c)\|_{\R_p^c}
 +\|(x_k^r)\|_{\R_p^r}\big\}
  \le c\,p'\, \|\sum_k x_k \|_p\,.\]
\end{theorem}

\pf Let  $(x_k)\in \R_p^d$. Then by \eqref{[1,2]},
 \[\|\sum_k x_k \|_p \le 2 \,
 \|(x_k)\|_{\R_p^d}\,.\]
To consider the second term on column norm, let $y_k=a_kD^{1/p}$
with $\E(a_k)=0$, and set $y=\sum_k y_k$. We deduce from
\cite[section 2]{ju-doob}(see also \cite[section 7]{jx-burk}) that
 \[\|y\|_p^2 = \|y^*y\|_{p/2}
 \le \|\E(y^*y)\|_{p/2}
 = \|\sum_{k}D^{1/p}\E(a_k^*a_k)D^{1/p}\|_{p/2} \,.\]
By density this implies that
 \[  \|\sum_k x_k \|_p \le
 \|(x_k)\|_{\R_p^c}\]
whenever $(x_k)\in\R_p^c$. Passing to adjoints, we get the same
inequality for the row subspace. Therefore, by triangle inequality
we find
 \[  \|\sum_k x_k \|_p\le 2\,
 \inf_{x_k=x_k^d+x_k^c+x_k^r}\,
 \big\{\|(x_k^d)\|_{\R_p^d} + \|(x_k^c)\|_{\R_p^c}
 +\|(x_k^r)\|_{\R_p^r}\big\}\,.\]
To prove the converse inequality we use duality. To this end note
that the infimum above is the norm of $(x_k)$ in the sum space
$\R_p^d+\R_p^c+\R_p^r$. By the duality between sums and
intersections, we have
 \[(\R_{p'}^d\cap\R_{p'}^c\cap\R_{p'}^r)^*=
 (\R_{p'}^d)^*+(\R_{p'}^c)^*+(\R_{p'}^r)^*\]
isometrically. However, by Lemma \ref{r-compl},
 \[(\R_{p'}^d)^*=\R_{p}^d\,,\quad (\R_{p'}^c)^*=\R_{p}^c\,,\quad
 (\R_{p'}^r)^*=\R_{p}^r\]
isomorphically. Therefore,
 \[(\R_{p'}^d\cap\R_{p'}^c\cap\R_{p'}^r)^*=
 \R_p^d+\R_p^c+\R_p^r\,.\]
Now let $ x_k\in L_p(\A_k)$ with $\E(x_k)=0$. Let $(y_k)\in
\R_{p'}^d\cap \R_{p'}^c\cap \R_{p'}^r$ such that
 \[\max\big\{\|(y_k)\|_{\R_{p'}^d}\,,\;
 \|(y_k)\|_{\R_{p'}^c}\,,\;\|(y_k)\|_{\R_{p'}^r}\big\}
 \le 1\,.\]
Then by Theorem \ref{ros},
 \[\|\sum_k y_k \|_{p'}\le c\,p'\,.\]
Thus, by orthogonality and the H\"older inequality
 \[\big|\sum_k \tr(y_k^*x_k)\big|
 =\big|\tr[(\sum_ky_k)^*(\sum_k x_k)]\big|
 \le c\,p'\,\big\|\sum_k x_k \big\|_p\,.\]
We then deduce the desired inequality. Hence the theorem is
proved. \qd

\smallskip

Now we give an application to random matrices. Recall that the
$e_{ij}$ denote the canonical matrix units of $B(\ell_2)$.

\begin{theorem}\label{randmat}
 Let $1<p<\8$ and $(x_{ij})$ be a finite matrix with
entries in $L_p(\M)$. Assume that the  $x_{ij}$ are independent
with respect to $\E$ and $\E(x_{ij})=0$. Then for $p\ge2$
 \begin{align*}
  &\big\|\sum_{ij}x_{ij}\ten e_{ij}
  \big\|_{L_p(\M\bar\ten B(\ell_2))}\,\sim_{cp}\\
  &~~\max\Big\{\big(\sum_{ij}\|x_{ij}\|_p^p\big)^{1/p}\,,\;
  \big(\sum_{j}\big\|\big[\sum_i\E(|x_{ij}|^2)\big]^{1/2}\big\|_p^p\big)^{1/p}\,,\;
   \big(\sum_{i}\big\|\big[\sum_j\E(|{x_{ij}}^*|^2)\big]^{1/2}\big\|_p^p\big)^{1/p}
    \Big\}
   \end{align*}
and for $p<2$
 \begin{align*}
  &\big\|\sum_{ij}x_{ij}\ten e_{ij}
  \big\|_{L_p(\M\bar\ten B(\ell_2))}\,\sim_{cp'}\\
  &~~\inf\Big\{\big(\sum_{ij}\|x^d_{ij}\|_p^p\big)^{1/p}+
  \big(\sum_{j}\big\|\big[\sum_i\E(|x^c_{ij}|^2)\big]^{1/2}\big\|_p^p\big)^{1/p}+
   \big(\sum_{i}\big\|\big[\sum_j\E(|{x^r_{ij}}^*|^2)\big]^{1/2}\big\|_p^p\big)^{1/p}
    \Big\},
   \end{align*}
where the infimum is taken over all decompositions
$x_{ij}=x^d_{ij}+x^c_{ij}+x^r_{ij}$ with mean zero elements
$x^d_{ij}$, $x^c_{ij}$ and $x^r_{ij}$, which, for each couple
$(i,j)$,  belong to the von Neumann subalgebra generated  by
$x_{ij}$.
\end{theorem}

\pf Assume that $(x_{ij})$ is an $n\times n$ matrix. Let ${\rm
Tr}$ be the usual trace on $B(\ell_2^n)$. Then $\f\ten{\rm Tr}$ is
a normal faithful positive functional  on $\M\bar\ten B(\ell_2^n)$
(which becomes a state if we wish by normalizing ${\rm Tr}$). The
conditional expectation from $\M\bar\ten B(\ell_2^n)$ onto
$\N\bar\ten B(\ell_2^n)$ is $\E\ten {\rm id}_{B(\ell_2^n)}$. It is
easy to see that $(x_{ij}\ten e_{ij})$ is independent with respect
to $\E\ten {\rm id}_{B(\ell_2^n)}$. Then the case $p\ge2$ follows
directly from Theorem \ref{ros}. Indeed, we have
 \be
  \big\|\sum_{ij} \E\ten{\rm id}_{B(\ell_2^n)}
  (|x_{ij}\ten  e_{ij}|^2)\big\|_{p/2}
  &=&
  \big\|\sum_{ij} \E(|x_{ij}|^2) \ten e_{jj}\big\|_{p/2}\\
  &=&\big(\sum_{j}
  \big\|\sum_i \E(|x_{ij}|^2)\big\|_{p/2}^{p/2}\big)^{2/p} \,.
  \ee
The same calculation applies to the second square function.

For the case $p<2$ we cannot formally apply Theorem \ref{dual}.
However, we can indeed follow the reduction argument of Theorem
\ref{dual} from Theorem \ref{ros}. For this let $(\A_{ij})$ be a
family of subalgebras independent over $\N$ such that $x_{ij}\in
L_p(\A_{ij})$. Accordingly, we define  $\tilde\R_p^c$ to be the
subspace of $\ell_p(L_p(\M,\E;\ell_2^c))$ consisting of $(y_{ij})$
such that $y_{ij}\in L_p(\A_{ij})$ and $\E(y_{ij})=0$. (Note that
$\ell_p$ and $\ell_2^c$ in $\ell_p(L_p(\M,\E;\ell_2^c))$ are in
$j$ and $i$, respectively; this corresponds to the second term in
the preceding maximum.) Then the proof of Lemma \ref{r-compl}
shows that $\tilde\R_p^c$ is complemented in
$\ell_p(L_p(\M,\E;\ell_2^c))$. Similarly, we introduce the
complemented diagonal and row subspaces $\tilde\R_p^d$ and
$\tilde\R_p^r$ of $\ell_p(\nz^2; L_p(\M))$ and
$\ell_p(L_p(\M,\E;\ell_2^r))$, respectively. The rest of the proof
is the same as that of Theorem \ref{dual}.\qd

\begin{rem}\label{random Sp}{\rm
 Applying Theorem \ref{randmat} to a Rademacher family
$(\eps_{ij})$ on a probability space $(\Om, \mu)$,  we get the
following well-known equivalence for $2\le p<\8$
 \begin{align*}
 &\big\|\sum_{ij}\eps_{ij}\,a_{ij}\,e_{ij}
  \big\|_{L_p(\Om; S_p)}\,\sim  \\
 &~~~\max\big\{\big(\sum_{j}\big(\sum_i|a_{ij}|^2
 \big)^{p/2}\big)^{1/p}\,,\;
 \big(\sum_{i}\big(\sum_j|a_{ij}|^2 \big)^{p/2}\big)^{1/p}
    \big\}
 \end{align*}
for all finite complex matrices $(a_{ij})$. Indeed, in this
special case the diagonal term
 $\big(\sum_{ij}|a_{ij}|^p\big)^{1/p}$
in the maximum is dominated by each of the two others (see
\eqref{p-2}). By duality, we get a similar equivalence for $1<p<2$
by replacing, as usual, the maximum by the corresponding infimum
(see \cite{lust-khin}). Note that $(\eps_{ij})$ can be replaced by
a standard Gaussian family.
 }\end{rem}

Applying the Rosenthal inequality to the independent sequences
contained in the examples of section \ref{independence}, we get
Khintchine type inequalities as in Corollaries \ref{scalar ros}
and \ref{khin-p2}. Because of their importance in applications, we
give some more details. For convenience, we group them together
into two remarks according to the tracial and non tracial cases.

\begin{rem}{\rm
 Let $\f$ be a normal faithful tracial state on $\M$, and let
$(x_k)$ be a sequence in $L_p(\M)$ such that
 \[\al_p=\inf_k\|x_k\|_p>0\quad\mbox{and}\quad
 \beta_p=\sup_k\|x_k\|_p<\8\]
for all $p<\8$. Assume that the $x_k$ are independent with respect
to $\f$ and $\f(x_k)=0$. Let $\A$  be another von Neumann algebra
and  $(a_k)\subset L_p(\A)$ a finite sequence. Then for
 $2\le p<\8$
 \[
 \big\|\sum_{k}a_k\ten x_k \big\|_{L_p(\A\bar\ten\M)}
 \sim
  \max\big\{\big\|(a_k)\big\|_{L_p(\A; \ell_2^c)}\,,\;
  \big\|(a_k)\big\|_{L_p(\A; \ell_2^r)}\big\}
 \]
  and for $1<p<2$
 \[
 \big\|\sum_{k}a_k\ten x_k \big\|_{L_p(\A\bar\ten\M)}
 \sim
  \inf\big\{\big\|(b_k)\big\|_{L_p(\A; \ell_2^c)}+
  \big\|(c_k)\big\|_{L_p(\A; \ell_2^r)}\big\},
 \]
where the infimum is taken over all decompositions $a_k=b_k+c_k$
in $L_p(\A)$. In both cases, the equivalence constants depend only
on $p$, $\al_p$ and $\beta_p$.

The first equivalence is a special case of Corollary
\ref{khin-p2}. The second then follows by duality. This statement
implies many known inequalities. For instance, if $(x_k)$ is a
Rademacher, Steinhauss or Gaussian sequence, we recover the
noncommutative Khintchine inequalities of Lust-Piquard/Pisier
\cite{LPP}. As far as for noncommutative independence, $(x_k)$ can
be a sequence of free Gaussians, $q$-Gaussians or free generators.
Then we get the corresponding inequalities already in
\cite{pis-ast} (except the $q$-case). It is worth to note that for
all these concrete examples, the second equivalence above holds
for $p=1$ too and the constant there is then controlled by a
universal one; moreover, in the noncommutative case (except
$q\neq-1$) the first equivalence is even true for $p=\8$ and the
constant is also universal (depending only on $q$ in the
$q$-case). We refer to \cite{pis-ast} for more information.
 }\end{rem}

\begin{rem}{\rm
 Here we consider only the quasi free CAR generators $(x_k)$
defined in  \eqref{aw}. Then
  for  $2\le p<\8$
 \begin{align*}
 &\big\|\sum_ka_k\ten D^{1/(2p)}x_kD^{1/(2p)}\big\|_p
 \sim\\
  &~~~\max\big\{\big\|\big(\sum_k(1-\mu_k)^{1/p}\mu_k^{1/p'}
  a_k^*a_k\big)^{1/2}\big\|_p\,,\;
  \big\|\big(\sum_k(1-\mu_k)^{1/p'}\mu_k^{1/p}
  a_ka_k^*\big)^{1/2}\big\|_p\big\}
 \end{align*}
 and for $1<p<2$
 \begin{align*}
 &\big\|\sum_{k}a_k\ten D^{1/(2p)}x_kD^{1/(2p)} \big\|_p
 \sim\\
  &~~~\inf\big\{\big\|\big(\sum_k(1-\mu_k)^{1/p}\mu_k^{1/p'}
  b_k^*b_k\big)^{1/2}\big\|_p+
  \big\|\big(\sum_k(1-\mu_k)^{1/p'}\mu_k^{1/p}
  c_kc_k^*\big)^{1/2}\big\|_p\big\}.
 \end{align*}
where the infimum is taken over all decompositions $a_k=b_k+c_k$
in $L_p(\A)$. Moreover, the equivalence constants depend only on
$p$.

This statement is a reformulation of \cite[Theorem 4.1]{xu-embed}.
Note that the case $p\ge2$ can be easily deduced from Corollary
\ref{khin-p2} and the other is again proved by duality. It is
shown in \cite{ju-araki} that the second equivalence remains true
for $p=1$. Let us point out that a similar statement holds for the
generalized circular system in \eqref{shlya}. In this case, all
constants are universal (see \cite{xu-gro}; see also \cite{jpx}
for the $q$-case). We should emphasize that all these Khintchine
type inequalities have interesting applications. In fact, they
play a crucial role in the recent works on the operator space
Grothendieck theorems and the complete embedding of Pisier's $OH$
into noncommutative $L_p$, see \cite{ju-OH, pisshlyak, xu-embed,
xu-gro}.
 }\end{rem}


\section{A variant using maximal functions}
 \label{maximal function variant}


We discuss in this section a version of the noncommutative
Rosenthal inequality where the diagonal norm of  $\ell_p(L_p(\M))$
is replaced by that of $L_p(\M;\ell_\infty)$. This is in perfect
analogy with the classical Burkholder inequality for commutative
martingales. Our argument is based on interpolation and the
resulting constant presents, unfortunately, a singularity as $p\to
2$. We need some facts on noncommutative $L_p(L_q)$. For our
purpose here we will need only the case where the second space
$L_q$ is $\ell_q$. The investigation of general noncommutative
$L_p(L_q)$ spaces will be pursued elsewhere.

\smallskip

Let us recall the definition of the spaces $L_p(\M;\ell_\infty)$
and $L_p(\M;\ell_1)$, $1\le p\le\8$. A sequence $(x_k)$ in
$L_p(\M)$ belongs to $L_p(\M;\ell_\8)$ iff $(x_k)$ admits a
factorization $x_k=ay_kb$ with $a, b\in L_{2p}(\M)$ and
$(y_k)\in\ell_\8(L_\8(\M))$. The norm of $(x_k)$ is then defined
as
 \begin{equation}\label{norm of Lplinfty}
  \|(x_k)\|_{L_p(\M;\ell_\8)}
 =\inf_{x_k=ay_kb}\, \|a\|_{2p}\,
 \|(y_k)\|_{\ell_\8(L_\8)}\, \|b\|_{2p}\,.
 \end{equation}
On the other hand,  $L_p(\M;\ell_1)$ is defined as the space of
all sequences $(x_k)\subset L_p(\M)$ for which there exist
$a_{kj}, b_{kj}\in L_{2p}(\M)$ such that
 \[x_k=\sum_ja_{kj}^*b_{kj}\,.\]
$L_p(\M;\ell_1)$ is equipped with the norm
 \[ \|(x_k)\|_{L_p(\M;\ell_1)}
 =\inf_{x_k=\sum_{j} a_{kj}^*b_{kj}}\,
 \big\|\sum_{k,j} a_{kj}^*a_{kj}\big\|_p^{1/2}\,
 \big\|\sum_{k,j} b_{kj}^*b_{kj}\big\|_p^{1/2} \,.\]
This norm has a description similar to that of
$L_p(\M;\ell_\infty)$:
 \begin{equation}\label{norm of Lpl1}
  \|x\|_{L_p(\M;\ell_1)}
 =\inf_{x_k=ay_kb}\, \|a\|_{2p}\,
 \|(y_k)\|_{L_\8(\M;\ell_1)}\, \|b\|_{2p}\,.
 \end{equation}
We refer to \cite{ju-doob} for more information (see also
\cite{jx-erg}). Now for $1<q<\8$ we define $L_p(\M;\ell_q)$ as a
complex interpolation space between $L_p(\M;\ell_\8)$ and
$L_p(\M;\ell_1)$:
 \[ L_p(\M;\ell_q)=
 [L_p(\M;\ell_\8),\; L_p(\M;\ell_1)]_{1/q} \,.\]
Our reference  for interpolation theory is \cite{bl}. The norm of
$L_p(\M;\ell_q)$ will be often denoted by $\|\,\|_{L_p(\ell_q)}$.
Let us note that if $\M$ is injective, this definition is a
special case of Pisier's vector-valued noncommutative $L_p$-space
theory \cite{pis-ast}. The following is also motivated by Pisier's
theory.

\begin{prop}\label{norm of Lplq}
 Let $(x_k)\subset L_p(\M)$. Then $(x_k)\in L_p(\M;\ell_q)$
iff $(x_k)$ admits a factorization $x_k=ay_kb$ with $a, b\in
L_{2p}(\M)$ and $(y_k)\in L_\8(\M; \ell_q)$. Moreover,
 \[ \|(x_k)\|_{L_p(\ell_q)}
 =\inf_{x_k=ay_kb}\, \|a\|_{2p}\,
 \|(y_k)\|_{L_\8(\ell_q)}\, \|b\|_{2p}\,.\]
\end{prop}

\pf Let $\tnorm{(x_k)}_{p,q}$ denote the infimum above. By
(\ref{norm of Lplinfty}) and (\ref{norm of Lpl1}), the trilinear
map $(a, (y_k), b)\mapsto (ay_kb)$ is contractive from
$L_{2p}(\M)\times L_\8(\M;\ell_q)\times L_{2p}(\M)$ to
$L_p(\M;\ell_q)$ for $q=\8$ and $q=1$, so is it for any $q\in(1,
\8)$ in virtue of interpolation. This yields
 \[\|(x_k)\|_{L_p(\ell_q)}\le\tnorm{(x_k)}_{p,q}\,.\]
To prove the converse we consider only the case where the state
$\f$ is tracial. The general case can be reduced to this one by
using Haagerup's reduction theorem as in \cite{xu-descrip}. Now
assume $\|x\|_{L_p(\ell_q)}<1$. Let $S=\{z\in\cz: 0\le{\rm Re}z\le
1\}$. Then there exists a sequence $(f_k)$ of continuous functions
from $S$ to $L_p(\M)$, analytic in the interior of $S$, such that
$f_k(1/q)=x_k$ and
 \[ \sup_{t\in\rz} \|(f_k(it))\|_{L_p(\ell_\8)}\le 1,
 \quad
 \sup_{t\in\rz} \|(f_k(1+it))\|_{L_p(\ell_1)}\le 1 . \]
By (\ref{norm of Lplinfty}) and (\ref{norm of Lpl1}), we have
factorizations
 \[ f_k(z)=a(z) y_k(z) b(z),
 \quad z\in \partial S \]
such that
 \[\|a(z)\|_{2p}\le 1,\quad \|b(z)\|_{2p}\le 1\]
and
 \[\|(y_k(it))\|_{L_\8(\ell_\8)}\le 1, \quad
 \|(y_k(1+it))\|_{L_\8(\ell_1)}\le 1. \]
Moreover,  we may assume that $a$, $b$ and $y$ are strongly
measurable on $\partial S$. Now fix $\eps>0$. Then by the
operator-valued Szeg\"o factorization \cite[Corollary
8.2]{px-survey}, we find two strongly measurable functions $\al,
\beta: S\to L_{2p}(\M)$, analytic in the interior, such that
  \[ \al(z)\al(z)^* = a(z)a(z)^*+ \eps \quad \mbox{and}
 \quad
 \beta(z)^*\beta(z) = b(z)^*b(z)+ \eps\,,
 \quad z\in\partial S\,. \]
Moreover, $\al(z)$ and $\beta(z)$ are invertible for every $z\in
S$. For $z\in\partial S$ let $u(z)$ and $v(z)$ be contractions in
$\M$ such that
 \[a(z)=\al(z)u(z)\quad\mbox{and}\quad b(z)=v(z)\beta(z)\,.\]
We then deduce
 \[f_k(z)=\al(z)u(z)y_k(z)v(z)\beta(z)\,.\]
Set $\tilde y_k(z)=u(z)y_k(z)v(z)$ for $z\in\partial S$. Since
$\al(z)$ and $\beta(z)$ are invertible, we have $\tilde
y_k(z)=\al(z)^{-1}f_k(z)\beta(z)^{-1}$. Thus $\tilde y_k$ is the
boundary value  of an analytic function in $S$, so $\tilde y_k$
itself may be viewed as an analytic function in $S$. Therefore, we
obtained an analytic factorization of $f_k$:
 \[f_k(z)=\al(z)\tilde y_k(z)\beta(z),\quad z\in S.\]
Moreover, we have the following estimates
 \[\|\al(z)\|_{2p}\le 1+\eps\,,\quad
 \|\beta(z)\|_{2p}\le 1+\eps\]
for any $z\in\partial S$ and
 \[\|(\tilde y_k(it))\|_{L_\8(\ell_\8)}\le 1, \quad
 \|(\tilde y_k(1+it))\|_{L_\8(\ell_1)}\le 1. \]
It then follows that
 \[\|\al(\frac1q)\|_{2p}\le 1+\eps\,,\quad
 \|\beta(\frac1q)\|_{2p}\le 1+ \eps\,\quad
 \|(\tilde y_k(\frac1q))\|_{L_\8(\ell_q)}\le 1.\]
Since $x_k=f_k(1/q)=\al(1/q)\tilde y_k(1/q)\beta(1/q)$, we deduce
 \[\|(x_k)\|_{L_p(\ell_q)}\le 1 +\eps\,.\]
Letting $\eps\to0$ yields $\|(x_k)\|_{L_p(\ell_q)}\le 1$. \qd

\begin{cor}
 Let $1\le p_0,p_1, q_0, q_1\le\8$ and $0<\tet<1$. Then
  \[ [L_{p_0}(\M; \ell_{q_0}),\; L_{p_1}(\M; \ell_{q_1})]_\tet
  =L_{p}(\M; \ell_{q})\]
 with equal norms, where $1/p=(1-\tet)/p_0+\tet/p_1$ and
$1/q=(1-\tet)/q_0+\tet/q_1$.
\end{cor}

\pf By Proposition \ref{norm of Lplq}, the trilinear map $(a,
(y_k), b)\mapsto (ay_kb)$ is contractive from $L_{2p_j}(\M)\times
L_\8(\M;\ell_{q_j})\times L_{2p_j}(\M)$ to
$L_{p_j}(\M;\ell_{q_j})$ for $j=0$ and $j=1$, so by interpolation
it is also contractive from $L_{2p}(\M)\times
L_\8(\M;\ell_{q})\times L_{2p}(\M)$ to $[L_{p_0}(\M;
\ell_{q_0}),\; L_{p_1}(\M; \ell_{q_1})]_\tet$. This, together with
Proposition \ref{norm of Lplq}, implies
 \[L_{p}(\M; \ell_{q})\subset [L_{p_0}(\M; \ell_{q_0}),\;
 L_{p_1}(\M; \ell_{q_1})]_\tet\,.\]
The converse inclusion is proved similarly as Proposition
\ref{norm of Lplq} by using the Szeg\"o factorization. We omit the
details. \qd

\begin{cor}\label{explicit norm of Lplq}
 Let $1\le p, q\le\8$.
 \begin{enumerate}[\rm(i)]
 \item $L_p(\M;\ell_p)=\ell_p(L_p(\M))$ isometrically.
 \item If $p\le q$,
 \[\|(x_k)\|_{L_p(\ell_q)}=\inf_{x_k=a y_kb}\,
 \|a\|_{2r}\,\|(y_k)\|_{\ell_q(L_q)}\,
 \|b\|_{2r}\]
for any $(x_k)\in L_p(\M;\ell_q)$, where $1/r=1/p-1/q$.
 \item If $p\ge q$,
  \[\|(x_k)\|_{L_p(\ell_q)}
  =\sup_{\|\al\|_{2s}\le1,\;\|\beta\|_{2s}\le1 }\,
 \|(\al x_k\beta)\|_{\ell_q(L_q)}\]
for any $(x_k)\in L_p(\M;\ell_q)$, where $1/s=1/q-1/p$.
 \end{enumerate}
\end{cor}

\pf (i) By definition the quality in question is true for $p=\8$
and $p=1$. For $1<p<\8$ we use the previous corollary to conclude
 \be
 L_p(\M;\ell_p)
 &=&[L_{\8}(\M; \ell_{\8}),\; L_{1}(\M; \ell_{1})]_{1/p}\\
 &=&[\ell_{\8}(L_{\8}(\M)),\; \ell_{1}(L_{1}(\M))]_{1/p}
 =\ell_p(L_p(\M))\,.
 \ee

(ii) Proposition \ref{norm of Lplq} may be  rewritten symbolically
as
 \[L_p(\M;\ell_q)=L_{2p}(\M)\,L_\8(\M;\ell_q)\,L_{2p}(\M)\,.\]
However, the H\"older inequality implies
 \[L_{2p}(\M)=L_{2r}(\M)\,L_{2q}(\M)=L_{2q}(\M)\,L_{2r}(\M)\,.\]
We thus deduce, by (i)
 \be
 L_p(\M;\ell_q)
 &=&L_{2r}(\M)\,L_{2q}(\M)\,L_\8(\M;\ell_q)\,L_{2q}(\M)\,L_{2r}(\M)\\
 &=&L_{2r}(\M)\,L_q(\M;\ell_q)\,L_{2r}(\M)
 =L_{2r}(\M)\,\ell_q(L_q(\M))\,L_{2r}(\M)\,;
 \ee
whence the desired result.

(iii) Given $(x_k)\in L_p(\M;\ell_q)$ we apply Proposition
\ref{norm of Lplq} to write $x_k=ay_kb$ with $a, b\in L_{2p}(\M)$
and $(y_k)\in L_\8(\M; \ell_q)$. Then for any $\al, \beta$ in the
unit ball of  $L_{2s}(\M)$, we have
 \begin{align*}
 \|(\al x_k\beta)\|_{\ell_q(L_q)}
 \le\|\al a\|_{2q}\,\|(y_k)\|_{L_\8(\ell_q)}\,\|b\beta\|_{2q}
 \le\|a\|_{2p}\,\|(y_k)\|_{L_\8(\ell_q)}\,\|b\|_{2p}\,.
 \end{align*}
Therefore,
 \[\sup_{\|\al\|_{2s}\le1,\;\|\beta\|_{2s}\le1 }\,
 \|(\al x_k\beta)\|_{\ell_q(L_q)}
 \le\|(x_k)\|_{L_p(\ell_q)}\,.\]
To prove the converse inequality, we use (ii) and duality. It
suffices to consider a finite sequence $(x_k)_{1\le k\le n}\subset
L_p(\M)$. Accordingly, we consider the $\ell_q^n$-valued
$L_p$-space $L_p(\M; \ell_q^n)$. We may also assume $p>q$. Then
 \[L_{p'}(\M; \ell_{1}^n)^*=L_p(\M; \ell_\8^n)\quad
 \mbox{and}\quad
 L_{p'}(\M; \ell_{\8}^n)^*=L_p(\M; \ell_1^n)\]
isometrically (see \cite{ju-OH} and \cite{jx-erg}). Using the
duality theorem on complex interpolation, we deduce
 \[L_{p'}(\M; \ell_{q'}^n)^*=L_p(\M; \ell_q^n).\]
Now let $(y_k)\in L_{p'}(\M; \ell_{q'}^n)$ be of norm less than
$1$. By (ii) we can write $y_k=az_kb$ with
 \[\|a\|_{2s}\le 1 ,\quad\|b\|_{2s}\le 1 ,\quad
 \|(z_k)\|_{\ell_{q'}(L_{q'})}\le 1 .\]
Then
 \[\big|\sum_k\tr(y_k^*x_k)\big|
 =\big|\sum_k\tr(z_k^*a^*x_kb^*)\big|
 \le\big\|(a^*x_kb^*)\big\|_{\ell_q(L_q)}\,;\]
whence the desired converse inequality. \qd

\begin{cor}\label{rowcol}
 Let $2\le p\le \infty$. Then
 \[ [L_p(\M;\ell_2^c),\; L_p(\M;\ell_2^r)]_{1/2} \subset
 L_p(\M;\ell_2) \,.\]
 \end{cor}

\pf Let $1/r=1/2-1/p$. We consider the map $T:(a,(x_k),b)\mapsto
(ax_kb)$. First, we note that
 \[ T: L_{\infty}(\M)\times L_p(\M; \ell_2^c)\times L_r(\M)\to
 \ell_2(L_2(\M)) \]
is a contraction because
 \be
 \sum_{k} \|ax_kb\|_2^2
 &\le&\|a\|_\8^2 \sum_{k} \tr(b^*x_k^*x_kb)=\|a\|_\8^2\,
 \tr\big((\sum_{k} x_k^*x_k)bb^*\big) \\
 &\le& \|a\|_\8^2\, \|\sum_{k} x_k^*x_k\|_{p/2}
 \|bb^*\|_{r/2} =\|a\|_\8^2\,
  \|(x_k)\|_{L_p(\M;\ell_2^c)}^2\, \|b\|_{r}^2 \, .
 \ee
Similarly, we see that
 \[ T: L_{r}(\M)\times L_p(\M; \ell_2^r)\times L_\infty(\M)\to
 \ell_2(L_2(\M)) \]
is a contraction. Thus by interpolation
 \[ T:L_{2r}(\M)\times
 [L_p(\M; \ell_2^c),\; L_p(\M; \ell_2^r)]_{1/2}\times
 L_{2r}(\M)\to \ell_2(L_2(\M)) \]
is a contraction. Then Corollary \ref{explicit norm of Lplq},
(iii) implies the assertion. \qd

\begin{rem}{\rm
 The  inclusion converse to that of Corollary \ref{rowcol}
holds too, so we have equality. This is a special case of the main
result from \cite{xu-descrip} (see also \cite{jupa-amalg} for more
general results of this type).
 }\end{rem}

Now we are ready to prove the version of the noncommutative
Rosenthal inequality in terms of maximal functions.

\begin{theorem}\label{ros2}
 Let $\N$ be a $\f$-invariant von Neumann subalgebra of $\M$ with
conditional expectation $\E$. Let $2< p<\infty$ and $(x_k)\subset
L_p(\M)$ be a sequence independent with respect to $\E$ such that
$\E(x_k)=0$. Then
 \begin{align*}
 &\|\sum_{k} x_k \|_p \le c_p
  \max\big\{\|(x_k)\|_{L_p(\ell_{\infty})}\,,\;
   \|(x_k)\|_{L_p(\M,\E;\ell_2^c)}\,,\;
   \|(x_k)\|_{L_p(\M,\E;\ell_2^r)}\big\} .
 \end{align*}
\end{theorem}

\pf If
 \[ \|(x_k)\|_{\ell_p(L_p)}
 < \max\big\{ \|(x_k)\|_{L_p(\M,\E;\ell_2^c)}\,,\;
   \|(x_k)\|_{L_p(\M,\E;\ell_2^r)}\big\},\]
then Theorem \ref{ros} implies
 \[ \|\sum_k x_k\|_p \le c\, p
 \max\big\{ \|(x_k)\|_{L_p(\M,\E;\ell_2^c)}\,,\;
   \|(x_k)\|_{L_p(\M,\E;\ell_2^r)}\big\} \,, \]
so we are done. It remains to consider the case where
 \[\max\big\{ \|(x_k)\|_{L_p(\M,\E;\ell_2^c)}\,,\;
   \|(x_k)\|_{L_p(\M,\E;\ell_2^r)}\big\}
   \le \|(x_k)\|_{\ell_p(L_p)}\,.\]
Again by Theorem \ref{ros}, we have
 \[  \|\sum_k x_k\|_p \le c\, p
 \|(x_k)\|_{\ell_p(L_p)} \,.\]
By the reiteration theorem, we deduce (with $\tet=2/p$)
  \[ L_p(\M;\ell_p)=[L_p(\M;\ell_{\infty}),\;L_p(\M;\ell_2)]_\theta \,.\]
This, together with Corollary \ref{explicit norm of Lplq} (i),
implies
 \[\big\|(x_k)\big\|_{\ell_p(L_p)}
 \le \big\|(x_k)\big\|_{L_p(\ell_{\infty})}^{1-\theta}\,
 \big\|(x_k)\big\|_{L_p(\ell_2)}^{\theta} \,.\]
Using Lemma \ref{unc} and (\ref{lower khintchine}), we have
 \begin{equation}\label{crp}
 \max\{\|(x_k)\|_{L_p(\M;\ell_2^c)}\,,\;
 \|(x_k)\|_{L_p(\M;\ell_2^r)} \}
 \le 2 \|\sum_k x_k\|_p \,.
 \end{equation}
Then by Corollary  \ref{rowcol}
 \[ \|(x_k)\|_{L_p(\ell_2)}\le 2 \|\sum_k x_k\|_p \,.\]
Combining these estimates we find  (after cancellation) that
 \[
 \|\sum_k x_k\|_p\le
 (c\, 2^{\theta}p)^{1/(1-\theta)}\, \|(x_k)\|_{L_p(\ell_\8)}\, .
 \]
 The theorem is thus proved with
  $c_p\le (c'p)^{p/(p-2)}$
 for $p>2$. In particular,  $c_p\le c''p$ for $p\ge4$. \qd

\smallskip

We take this opportunity to present the same improvement in the
context of the noncommutative Burkholder inequality of
\cite{jx-burk}. Namely, we want to replace the norm
$\|(dx)\|_{\ell_p(L_p)}$ in the following inequality by
$\|(dx)\|_{L_p(\ell_\8)} $:
 \[\|x\|_p\le
 c_p\max\big\{\|(dx)\|_{\ell_p(L_p)}\,,\;\|x\|_{h_p^c}\,,\;
 \|x\|_{h_p^r}\big\}\]
for any  noncommutative martingale $x=(x_k)$  with respect to an
increasing filtration $(\E_k)$ of normal faithful conditional
expectations. Here  $dx=(dx_k)$ denotes the difference sequence of
$x$ and
 \[\|x\|_{h_p^c}
 =\big\|\big(\sum_k\E_{k-1}(|dx_k|^2)\big)^{1/2}\|_{p}\,,\quad
 \|x\|_{h_p^r}=\|x^*\|_{h_p^c}\,.\]
We refer to \cite{jx-burk} for more details. Note that $c_p\le
c\,p$ according to \cite{ran-weak}, which improves the original
estimate $c_p\le c\,p^2$ from \cite{jx-burk}.

\begin{theorem} Let $2<p<\8$. Then for any noncommutative bounded
$L_p$-martingale $x$ we have
 \[\|x\|_p\le
 c'_p\max\big\{\|(dx)\|_{L_p(\ell_\8)}\,,\;\|x\|_{h_p^c}\,,\;
 \|x\|_{h_p^r}\big\}.\]
\end{theorem}

\pf This proof is almost the same as that of the previous theorem.
The only difference is that the martingale analogue of (\ref{crp})
is now obtained by using the lower estimate in the noncommutative
Burkholder-Gundy inequality (see \cite{jx-const} for the optimal
order of the constant):
 \[\max\{\|(dx)\|_{L_p(\M;\ell_2^c)}\,,\;
 \|(dx)\|_{L_p(\M;\ell_2^r)} \}
 \le c\,p\, \|x\|_p \,.\]
We omit the details. The resulting order of the constant $c'_p$ is
the same as that of $c_p$ in the previous theorem. \qd

\begin{rem}{\rm
We can also improve the lower estimates in the noncommutative
Burkholder/Rosenthal inequalities for $1<p<2$, by replacing the
diagonal term $\ell_p(L_p)$ by $L_p(\ell_1)$. For instance, under
the assumptions of Theorem \ref{dual} we have
 \begin{align*}
 \inf_{x_k=x_k^d+x_k^c+x_k^r}\,
 \big\{\|(x_k^d)\|_{\tilde\R_p^d} + \|(x_k^c)\|_{\R_p^c}
 +\|(x_k^r)\|_{\R_p^r}\big\}
  \le c_p \|\sum_k x_k \|_p \, ,
 \end{align*}
where $\tilde\R_p^d$ is the subspace of $L_p(\M;\ell_1)$
consisting of all $(x_k)$ such that $x_k\in L_p(\A_k)$ with
$\E(x_k)=0$. The proof is similar to that of Theorem \ref{dual}
via duality. The complementation of the space $\tilde\R_p^d$
follows from the noncommutative  Doob inequality in
\cite{ju-doob}.
 }\end{rem}


\section{The nonfaithful case}\label{nonfaithful}


Nonfaithful filtrations of von Neumann subalgebras, so nonfaithful
conditional expectations,  occur very naturally in operator
algebra theory. The simplest example is the natural filtration
$(\Ma_n)_{n\ge1}$ of $B(\ell_2)$ given by the algebras $\Ma_n$ of
matrices $(a_{ij})$ such that $a_{ij}=0$ if $\max(i, j)>n$. On the
other hand, the notion of nonfaithful copies in a tensor product
of von Neumann algebras is important in the context of iterated
ultraproducts of von Neumann algebras.

The aim of this section is to extend Theorem \ref{ros} to the case
of nonfaithful conditional expectations. We start with the
relevant notion. $\M$ is still assumed  $\si$-finite and equipped
with a normal faithful state $\f$.  Let $\N$ be a w*-closed
involutive (not necessarily unital) subalgebra of $\M$. Let  $e$
be the unit of $\N$, so $e$ is a projection of $\M$. Again, assume
that $\N$ is $\f$-invariant (i.e., $\si_t^\f(\N)\subset\N$ for all
$t\in\rz$). With these assumptions we still have  a normal
conditional expectation $\E_\N:\M\to \N$ with support equal to $e$
such that $\f\circ\E_\N=\f_e$, where $\f_e=e\f e$. Like in the
faithful case, $\E_\N$ extends to a contractive projection from
$L_p(\M)$ onto $L_p(\N)$ for every $p\ge1$. We refer to
\cite{jx-burk} for more details.

Now, we consider a sequence $(\A_k)$ of $\f$-invariant w*-closed
involutive subalgebras of $\M$ containing $\N$. Let us denote by
$r_k$ the unit of $\A_k$. We will say that the algebras $\A_k$ are
\emph{independent over $\N$} or \emph{with respect to $\E_\N$} if
 \begin{enumerate}[(i)]
 \item  the projections $s_k=r_k-e$ are mutually orthogonal;
 \item  for every $k$, $\E_\N(xy)= \E_\N(x)\E_\N(y)$ holds for all
$x\in \A_k$ and $y$ in the w*-closed involutive subalgebra
generated by $(\A_j)_{j\neq k}$.
 \end{enumerate}
Note that in this case $(e\A_k e)$ is faithfully independent over
$\N$ in the sense of  section \ref{independence}. A sequence
$(x_k)\subset L_p(\M)$ is called \emph{independent with respect to
$\E_\N$} if there exists a sequence $(\A_k)$ of subalgebras
independent with respect to $\E_\N$ such that $x_k\in L_p(\A_k)$.

The new ingredient for the nonfaithful version of the
noncommutative Rosenthal inequality is a separate treatment of the
corners. In the rest of this section we will assume that $(\A_k)$
is independent with respect to $\E=\E_\N$ and keep the preceding
notations.

\begin{lemma}\label{corn} Let $2\le p<\infty$ and $x_k\in
L_p(\A_k)$.   Then
 \[ \|\sum_k s_kx_ke\|_p
 \le c\,\sqrt p\, \max \big\{\|(s_kx_ke)\|_{\ell_p(L_p)},\;
  \|(s_kx_k)\|_{L_p(\M,\E;\ell_2^c)}\big\} .
  \]
\end{lemma}

\pf Let $x=\sum_k s_kx_ke$. By the orthogonality of the $s_k$, we
obtain
 \[\|x\|_p^2 = \|\sum_k ex_k^*s_kx_ke \|_{p/2}
 \le\|\sum_k \E(x_k^*s_kx_k) \|_{p/2}
 + \|\sum_k ex_k^*s_kx_ke-\E(x_k^*s_kx_k) \|_{p/2}
 \,. \]
Note that $y_k=ex_k^*s_kx_ke-\E(x_k^*s_kx_k) \in e\A_ke$ and
satisfies $\E(y_k)=0$. As observed before, the sequence $(e\A_ke)$
is faithfully independent over $\N$. Now, we follow the proof of
Theorem \ref{ros}. If $2\le p\le 4$, we deduce from \eqref{[1,2]}
that
 \begin{align*}
 \|\sum_k y_k\|_{p/2}
 &\le 2 \ez \,\|\sum_k \eps_ky_k\|_{p/2}
 \le 2 \big(\sum_{k} \|y_k\|_{p/2}^{p/2}\big)^{2/p}
 \le 4\big(\sum_{k} \|s_kx_ke\|_{p}^{p}
 \big)^{2/p}
 \,.
 \end{align*}
For $4<p<\infty$, we deduce from Theorem \ref{ros} applied to
$(y_k)\subset L_{q}(e\M e)$ with $q=p/2$ and Lemma \cite[Lemma
5.2]{jx-burk} that
 \be
 \|\sum_k y_k\|_{q}
 &\le&  cp
 \max\Big\{ \big\|(y_k)\big\|_{\ell_q(L_q)}\,,\;
  \big\|\big(\sum_k \E(y_k^*y_k)\big)^{1/2}\big\|_{q}\Big\} \\
 &\le& cp
 \max\Big\{\big\|(s_kx_ke)\big\|^2_{\ell_p(L_p)}\,,\;
  \big\|\sum_k \E(|s_kx_ke|^4)\big\|_{p/4}^{1/2}\Big\} \\
 &\le&   cp
 \max\Big\{\big\|(s_kx_ke)\big\|^2_{\ell_p(L_p)}\,,
 \big(\sum_k \|s_kx_ke\|_p^p\big)^{1/(p-2)}\,
  \big\|\sum_k \E(|s_kx_ke|^2)\big\|_{p/2}^{(p/2-2)/(p-2)}
  \Big\}.
  \ee
Then the assertion follows by homogeneity.
 \qd

\smallskip

The nonfaithful version of the Rosenthal inequality for $p\ge 2$
has the same form as Theorem \ref{ros}.

\begin{theorem}\label{fp} Let $2\le p<\infty$ and $x_k\in
L_p(A_k)$. Set $y_k= x_k-\E(x_k)$. Then
 \begin{align*}
  \|\sum_k x_k \|_p\sim_{cp}
  \max\big\{\|\sum_k \E(x_k)\|_p,\;\|(y_k)\|_{\ell_p(L_p)},\;
  \|(y_k)\|_{L_p(\M,\E;\ell_2^c)},\;
  \|(y_k)\|_{L_p(\M,\E;\ell_2^r)}\big\}\,.
  \end{align*}
\end{theorem} \allowdisplaybreaks

\pf Since
 \[\|\sum_k x_k \|_p\le \|\sum_k \E(x_k)\|_p
 + \|\sum_k y_k\|_p\,,\]
we need only to estimate the second term on the right. Since $y_k$
is supported by $r_k$ and $s_k=r_k-e$ for each $k$,  we have
 \[\|\sum_k y_k\|_p\le \|\sum_k s_ky_ks_k\|_p
 +\|\sum_k s_ky_ke\|_p +\|\sum_k ey_ks_k\|_p
 +\|\sum_k ey_ke\|_p\,.\]
By the mutual orthogonality of the $s_k$,
 \[\|\sum_k s_ky_ks_k\|_p
 =\big(\sum_k \|s_ky_ks_k\|_p^p\big)^{1/p}
 \le \|(y_k)\|_{\ell_p(L_p)}\,.\]
On the other hand, by Lemma \ref{corn},
 \be
 \|\sum_k s_ky_ke\|_p
 &\le& c\sqrt p\,
 \max \big\{\|(s_ky_ke)\|_{\ell_p(L_p)},\;
  \|(s_ky_k)\|_{L_p(\M,\E;\ell_2^c)}\big\}\\
  &\le& c\sqrt p\,
 \max \big\{\|(y_k)\|_{\ell_p(L_p)},\;
  \|(y_k)\|_{L_p(\M,\E;\ell_2^c)}\big\}.
 \ee
Passing to adjoints, we get the same estimate for another term on
the corners. To deal with the last term, we recall that the
algebras $e\A_ke$ are faithfully independent over $\N$. Thus
Theorem \eqref{ros} applies to $(ey_ke)$:
 \be
 \|\sum_k ey_ke\|_p
 &\le& c p
 \max\big\{\|(ey_ke)\|_{\ell_p(L_p)},\;
  \|(ey_ke)\|_{L_p(\M,\E;\ell_2^c)},\;
  \|(ey_ke)\|_{L_p(\M,\E;\ell_2^r)}\big\}\\
  &\le& c p
 \max\big\{\|(y_k)\|_{\ell_p(L_p)},\;
  \|(y_k)\|_{L_p(\M,\E;\ell_2^c)},\;
  \|(y_k)\|_{L_p(\M,\E;\ell_2^r)}\big\}.
 \ee
Combing the preceding inequalities, we obtain the upper estimate.
The lower estimate is proved in the same way as in the faithful
case. \qd

\begin{exam}\label{corner} {\rm
 Nonfaithful independence occurs
naturally in the context of conditional expectations with respect
to corners. Let $\M$ be a von Neumann algebra, $e$ a projection
and $(r_k)$  a family of  projections such that $e\le r_j$ and
such that the $s_k$ are mutually orthogonal, where $s_k=r_k-e$.
Consider
 \[ \N =e\M e \quad\mbox{and} \quad \A_k = r_k\M r_k \,.\]
The conditional expectation associated with $\N$ is given by
$\E(x)=exe$. Then the $\A_k$ are independent over $\N$. This
situation occurs for example on a tensor product $\M=\B^{\ten_n}$,
where $e=f_1\ten \cdots \ten f_n$ and $r_k= f_1\ten \cdots
f_{k-1}\ten 1\ten f_{k+1}\ten \cdots \ten f_n$ with $f_k$
projections in $\B$.
 }\end{exam}

\begin{rem}{\rm
There exists, of course, a nonfaithful version of the Rosenthal
inequality for $1<p\le 2$. We keep the same assumptions as before.
The main technical difference is that we have to introduce two
extra spaces
 \[ R_p(se)=\big\{ \sum_k s_kx_ke\;:\; x_k\in L_p(\A_k),\;
 \E(x_k)=0\big\}\]
and
  \[R_p(es)  = \big\{ \sum_k ex_ks_k\;:\; x_k\in L_p(\A_k),\;
  \E(x_k)=0 \big\}\,.
  \]
It is easy to show that they are complemented in
 \[ \big\{\sum_k x_k \;:\; x_k \in L_p(\A_k),\;
  \E(x_k)=0\big\} \,.\]
Thanks to Lemma \ref{corn}, we are able to describe the dual
$R_{p'}(se)$ of $R_p(se)$ as an intersection of two terms, an
$\ell_{p'}$-term and a column square function. Using the duality
argument from the proof of Theorem \ref{dual}, we deduce
 \[\big\|\sum_k s_kx_ke\big\|_{p} \,\sim_{c\sqrt{p'}}
  \inf_{s_kx_ke=s_kx_k^de+s_kx_k^ce}\,
  \big\|(s_kx_k^de)\big\|_{\ell_p(L_p)} +
   \big\|(s_kx_k^c)\big\|_{L_p(\M,\E;\ell_2^c)}\,.\]
A similar result holds for $R_p(es)$. Now let $x_k\in L_p(A_k)$
with $\E(x_k)=0$. Then
 \[\big\|\sum_k x_k\big\|_p\, \sim_{c}\,
  \max\big\{\big\|(s_kx_ks_k)\big\|_{\ell_p(L_p)},\;
  \big\|\sum_k s_kx_ke\big\|_p,\;
  \big\|\sum_k ex_ks_k\big\|_p,\;
  \big\|\sum_k ex_ke\big\|_{p}
 \big\}.\]
The second and third terms were already treated. However,  the
last term  is the faithful part, so can be dealt with according to
Theorem \ref{dual}, which yields an equivalence with an infimum.
This complicated expression involving maximum and infimum is
particularly interesting in connection with independent copies (as
in \cite{ju-araki}). In this case, the expressions are symmetric.
This formula can be used to prove that subsymmetric sequences in
$L_p(\M)$, $1<p\le 2$, are symmetric (see \cite{juros} for more
details).
 }\end{rem}


\section{Symmetric subspaces of noncommutative  $L_p$}
 \label{symmetric subspaces}


In this section, we present some applications of the
Burkholder/Rosenthal inequalities to the study of symmetric
subspaces of noncommutative $L_p$-spaces both in the category of
Banach spaces and in that of operator spaces. The results obtained
are the noncommutative or operator space analogues of the
corresponding results in \cite{jmst}. Thus we will follow
arguments in \cite{jmst} in many cases. It will be convenient to
state these results in parallel for both categories, which will
also ease comparing and understanding them. All unexplained Banach
space terminologies used in the sequel can be found in
\cite{LT-I}. We refer to \cite{er-book, pis-intro} for background
on operator spaces and completely bounded maps and to
\cite{pis-ast, jnrx-OLp} for the operator space structure of
noncommutative $L_p$-spaces. In this paper we will focus on
subspaces of these spaces. In this situation we will only need the
following fact from \cite{pis-ast}: If $X$ and $Y$ are subspaces
of $L_p(\M)$, $1\le p\le\8$, then the cb-norm of a linear map
$T:X\to Y$ is given by
 \[ \|T\|_{cb}=\|{\rm id}_{S_p}\ten T:S_p(X)\to
 S_p(Y)\| \,.\]
Here $S_p(X)$ denotes the closure of $S_p\ten X$ in
$L_p(B(\ell_2)\bar\ten \M)$. In other words, the cb-norm is
calculated with matrix-valued coefficients instead of
scalar-valued coefficients for the usual norm $\|T\|$. It is then
straightforward to transfer to this setting all Banach space
terminologies concerning bases, basic sequences, {\it etc.}. For
instance, a basic sequence $(x_k)\subset X$ is said to be
completely unconditional if there exists a constant $\la$ such
that
 \[\|\sum_k\eps_ka_k\ten x_k\|\le\la\,\|\sum_ka_k\ten x_k\|\]
for all $a_k\in S_p$ and $\eps_k=\pm1$. Similarly, a FDD (finite
dimensional decomposition) $(F_k)$ of $X$ is said to be completely
unconditional  if  there exists a constant $\la$ such that
 \[\|\sum_k\eps_ka_k\ten x_k\|\le\la\,\|\sum_ka_k\ten x_k\|\]
for all $x_k\in F_k$,  $a_k\in S_p$ and $\eps_k=\pm1$.

\smallskip

The von Neumann algebras considered in this section and the next
one may  be non $\si$-finite. However, since we will often
consider sequences or separable subspaces in $L_p(\M)$, it is easy
to bring $\M$ to a $\si$-finite subalgebra  (see also Remark
\ref{general VN}).

\begin{lemma}\label{factor}
 Let $\M$ be a hyperfinite type {\rm III}$_{\la}$ factor
with $0\le \la\le 1$. Let $1<p<\infty$. Then $L_p(\M)$ has a
completely unconditional FDD.
\end{lemma}

\pf  In the range $0<\la\le 1$, we may assume that $\M$ is an
ITPFI factor. In general (including $\la=0$), we can always find a
normal faithful state $\f$, and an increasing sequence of finite
dimensional $\f$-invariant subalgebras $\M_n$ with conditional
expectations $\E_n: \M\to \M_n$ (see \cite{jrx-OLp2}). This yields
a martingale structure on $\M$. We define the difference operators
$\D_n=\E_n-\E_{n-1}$ where $\E_0=0$. Note that the spaces
$F_n=\D_n(L_p(\M))$ are finite dimensional and every element can
be written uniquely as $x=\sum_n \D_n(x)$. Thus $L_p(\M)$ has a
FDD. The complete unconditionality of this decomposition  means
that all maps $T_{\eps}=\sum_n \eps_n \D_n$ are completely bounded
uniformly in $\eps_n=\pm1$. Namely, the maps ${\rm id}_{S_p}\ten
T_{\eps}$ are uniformly bounded, i.e., there exists a constant $c$
such that
  \[
  \big\|\sum_n \eps_n ({\rm id}_{S_p}\ten \D_n)(x)\big\|_p \le  c\,
  \big\|\sum_n ({\rm id}_{S_p}\ten \D_n)(x)\big\|_p \,.\]
holds for all choices of signs $(\eps_n)$ and $x\in
L_p(B(\ell_2)\bar\ten \M)$. But this inequality  is a direct
consequence of the noncommutative Burkholder-Gundy inequalities
\cite{px-BG, jx-burk}. Moreover, the constant $c$ depends only on
$p$. \qd

\begin{theorem}\label{hyp}
 Let $\M$ be a hyperfinite von Neumann algebra. Let  $2<p<\infty$,
and let $(x_n) \subset L_p(\M)$ be a sequence of unit vectors,
which converges weakly to $0$. Then there exist constants $0\le
\al,\beta\le 1$, depending only on $(x_n)$, and a subsequence
$(\tilde x_n)$ of $(x_n)$ such that
 \[
 \big\|\sum_na_n\ten\tilde x_n\big\|_p\sim_{c_p}
 \max\big\{\big(\sum_n\|a_n\|_p^p\big)^{1/p}\,,\;
 \al\big\|\big(\sum_na_n^*a_n\big)^{1/2}\big\|_p\,,\;
 \beta\big\|\big(\sum_na_na_n^*\big)^{1/2}\big\|_p\big\}\]
 holds for all finite sequences $(a_n)\subset S_p$.
 \end{theorem}

\pf The first part of the proof is to show that we can reduce our
problem to the case where $L_p(\M)$ has a completely unconditional
FDD. To this end we first use a standard procedure to reduce $\M$
to a von Neumann algebra with separable predual (see
\cite[Appendix]{ggms}). Indeed, assume that $\M$ is $\si$-finite
and let $\f$ be a normal faithful state on $\M$. Let $A\subset \M$
be a countable subset, and let $\M_A$ be the von Neumann
subalgebra generated by $\si_t^{\f}(a)$ with $a\in A$ and $t\in
\qz$. Then $\M_A$  has separable predual. Moreover, $\M_A$ is
$\f$-invariant. Consequently, there is a normal faithful
conditional expectation from $\M$ onto $\M_A$, thus $L_p(\M_A)$ is
a complemented subspace of $L_p(\M)$. Now writing each $x_n$ as a
convergent series of elements from $\M D^{1/p}$:
 $x_n=\sum_{k}a_{nk}D^{1/p}\,$,
we can take $\{a_{nk}: n, k\in\nz\}$ as $A$. Then  $x_n\in
L_p(\M_A)$. Therefore, replacing $\M$ by $\M_A$,  we may assume
$\M_*$ separable.

Now if $\M$ is semifinite, then by \cite[Theorem~3.4]{pis-ast}
$\M$ has an increasing filtration of finite dimensional
subalgebras; so as in the proof of Lemma \ref{factor}, we deduce
that $L_p(\M)$ has a completely unconditional FDD. To treat the
case where $\M$ is of type III, we use another standard trick in
order to ensure that we may work with  a factor. To this end, we
consider the crossed product $\R=\overline{\ten}_{\nen} (\M,\f)
\rtimes G$ between the infinite tensor product
$\overline{\ten}_{\nen} (\M,\f)$ and the discrete group $G$ of all
finite permutations on $\nz$. Any finite permutation acts on the
infinite tensor product by shuffling the corresponding
coordinates. Clearly, we also have a normal faithful conditional
expectation $\E :\R\to \M$ obtained by first projecting onto the
identity element of $G$ and then to the first component in the
infinite tensor product. This implies that $L_p(\M)$ can be
identified as a (complemented) subspace of $L_p(\R)$. On the other
hand, according to \cite[Proof of Theorem 2.6]{haag-win}, $\R$ is
a hyperfinite factor. Thus $\R$ is of type III$_{\la}$ for some
$0\le \la \le 1$ (see \cite{con-classi, haag-bicent}). Therefore,
Lemma \ref{factor} implies that $L_p(\R)$ has a completely
unconditional FDD given by a filtration $(\E_k)$ of normal
faithful conditional expectations. In the remainder of the proof,
replacing $\M$ by $\R$ if necessary, we may assume that $L_p(\M)$
itself has this FDD.

The second part of the proof follows very closely its commutative
model (see \cite[Theorem 1.14]{jmst}). Using the gliding hump
procedure, we may find a perturbation of a subsequence
$(\hat{x}_n)$ and a corresponding subsequence $(\hat{\E}_k)$ such
that

\begin{samepage} \begin{enumerate}[(i)]
 \item  $\hat{\E}_n(\hat{x}_n)=\hat{x}_n$;
 \item  $\hat{\E}_{n}(\hat{x}_k)=0$ for all $k>n$;
 \item  $\lim_k \hat{\E}_n(\hat{x}_k^*\hat{x}_k)=y_n$ and
 $\|\hat{\E}_n(\hat{x}_k^*\hat{x}_k)-y_n\|_{p/2}
 \le\eps 2^{-k}$ for $k>n$;
 \item $\lim_k \hat{\E}_n(\hat{x}_k\hat{x}_k^*)=z_n$ and
 $\|\hat{\E}_n(\hat{x}_k\hat{x}_k^*)-z_n\|_{p/2}
 \le \eps 2^{-k}$ for $k>n$.
 \end{enumerate}\end{samepage}
Here $\eps>0$ is arbitrarily given and will be chosen after
knowing the $y_n$'s. It follows immediately from (iii) that
$(y_n)$ is a bounded $L_{p/2}$-martingale with respect to
$(\hat\E_n)$. Since $p/2>1$,  $(y_n)$ converges to some  $y\in
L_{p/2}(\M)$. Similarly, we obtain that $(z_n)$ converges to some
$z\in L_{p/2}(\M)$.  We define $\al=\|y\|_{p/2}^{1/2}$ and
$\beta=\|z\|_{p/2}^{1/2}$. Passing to subsequences of
$(\hat{x}_n)$ and  $(\hat{\E}_k)$ if necessary, we may further
assume
 \[ \big\|\hat{\E}_{n-1}(\hat{x}_{n}^*\hat{x}_{n})- y\big\|_{p/2}
  \le 2^{-(n+1)} \|y\|_{p/2}\,,  \quad
   \big\|\hat{\E}_{n-1}(\hat{x}_{n}\hat{x}_{n}^*)- z\big\|_{p/2}
  \le 2^{-(n+1)} \|z\|_{p/2}  \,. \]
Note that (i) and (ii) imply that $(\hat x_n)$ is a martingale
difference sequence with respect to $(\hat\E_n)$. Thus applying
the noncommutative Burkholder inequality \cite{jx-burk}, we find,
for any $a_n\in S_p$,
 \begin{align*}
  \big\|\sum_n a_n \ten \hat{x}_n \big\|_p
  \sim_{c_p} &\big(\sum_n \big\|a_n\ten
  \hat{x}_n\big\|_p^p\big)^{1/p}
   + \big\|\sum_n a_n^*a_n \ten
  \hat{\E}_{n-1}(\hat{x}_{n}^*\hat{x}_{n})\big\|_{p/2}^{1/2}\\
  &~~   + \big\|\sum_n a_na_n^* \ten
  \hat{\E}_{n-1}(\hat{x}_{n}\hat{x}_{n}^*)\big\|_{p/2}^{1/2} \,.
  \end{align*}
From perturbation, we have $1/2 \le \|\hat{x}_n\|_p\le 2$, so the
first diagonal term on the right is fine. On the other hand, the
triangle inequality implies
 \begin{align*}
 &\big\|\sum_n a_n^*a_n \ten
  \hat{\E}_{n-1}(\hat{x}_{n}^*\hat{x}_{n}) - \sum_n a_n^*a_n
  \ten  y \big\|_{p/2}\\
  &\quad  \le  \sum_n \big\|a_n^*a_n\big\|\,
  \big\|\hat{\E}_{n-1}(\hat{x}_{n}^*\hat{x}_{n})
  -y\big\|_{p/2} \\
  &\quad  \le  \frac12  \|y\|_{p/2}\, \sup_n \big\|a_n^*a_n\big\|_{p/2}
  \le\frac{\al^2}2\, \big\|\sum_n a_n^*a_n\big\|_{p/2} \,.
  \end{align*}
It follows that
 \[\Big|\big\|\sum_n a_n^*a_n \ten
  \hat{\E}_{n-1}(\hat{x}_{n}^*\hat{x}_{n})\big\|_{p/2}-
   \big\|\sum_n a_n^*a_n\ten y\big\|_{p/2}\Big|\le
   \frac{\al^2}2\, \big\|\sum_n a_n^*a_n\big\|_{p/2}\,.\]
However,
 \[\big\|\sum_n a_n^*a_n\ten y\big\|_{p/2}=
 \al^2\,\big\|\sum_n a_n^*a_n\big\|_{p/2}\,.\]
Therefore, we deduce
  \[\big|\big\|\sum_n a_n^*a_n \ten
  \hat{\E}_{n-1}(\hat{x}_{n}^*\hat{x}_{n})\big\|_{p/2}\sim_c
  \al^2\,\big\|\sum_n a_n^*a_n\big\|_{p/2} \,.\]
The same argument applies to the last term on the row norm.
Keeping in mind that $(\hat x_n)$ is a perturbation of a
subsequence $(\tilde x_n)$ of $(x_n)$ and going back to this
subsequence, we get the announced result. \qd

\smallskip

As a first application, we present an operator space version of
the Kadec-Pe\l zsy\'{n}ski alternative. For this we need some
notation from the theory of operator spaces. The spaces $C_p$ and
$R_p$ are defined as the column and row subspaces of $S_p$,
respectively. Namely,
 \[ C_p =\overline{{\rm span}}\,\{ e_{k1}\; :\; k\in \nz\}
 \quad \mbox{and} \quad R_p=\overline{{\rm span}}\,\{ e_{1k}\; :\; k\in
 \nz\}\,.\]
Note that as Banach spaces, $C_p$ and $R_p$ are isometric to
$\ell_2$ by identifying both $(e_{k1})$ and $(e_{1k})$ with the
canonical basis $(e_k)$ of $\ell_2$. We will adopt this
identification in the sequel. This permits us to consider the
intersection $C_p\cap R_p$. Recall that the operator space
structures of these spaces are determined as follows. For any
finite sequence $(a_k)\subset S_p$,
 \[\|\sum_k a_k\ten e_k\|_{S_p(C_p)}
 =\|(\sum_ka_k^*a_k)^{1/2}\|_p\,,\quad
 \|\sum_k a_k\ten e_k\|_{S_p(R_p)}
 =\|(\sum_ka_ka_k^*)^{1/2}\|_p \]
and
  \[\|\sum_k a_k\ten e_k\|_{S_p(C_p\cap R_p)}
 =\max\big\{\|(\sum_ka_k^*a_k)^{1/2}\|_p\,,\;
 \|(\sum_ka_ka_k^*)^{1/2}\|_p\big\}.\]

Recall that a sequence $(x_k)$ in a Banach space $X$ is said to be
semi-normalized if $\inf_k\|x_k\|>0$ and $\sup_k\|x_k\|<\8$.

\begin{cor}\label{ssp} Assume that $\M$ is hyperfinite.
Let $2\le p<\infty$ and $(x_n)\subset L_p(\M)$ be a
semi-normalized sequence which converges to $0$ weakly. Then
$(x_n)$ contains a subsequence $(\tilde x_n)$ which is completely
equivalent to the canonical basis of $\ell_p$, $C_p$, $R_p$ or
$C_p\cap R_p$.
\end{cor}

\pf Assume $p>2$. Let  $(\tilde x_n)$ be the subsequence from
Theorem \ref{hyp}. If   $\al=\beta=0$, then $(\tilde x_n)$ is
completely equivalent to the basis of $\ell_p$. If $\al>0$ and
$\beta=0$, then we find a copy of $C_p$ by virtue of (\ref{p-2}).
Similarly, if $\al=0$ and $\beta>0$ it turns out to be $R_p$. The
case $\al>0$ and $\beta>0$ yields $C_p\cap R_p$. \qd

\smallskip

A basis $(x_k)$ of $X\subset L_p(\M)$ is called symmetric if there
exists a positive constant $\la$ such that
  \[
  \|\sum_k \eps_k \al_{\pi(k)} x_k \|_p \le\la\,
  \|\sum_k  \al_{k } x_k \|_p
  \]
holds for finite sequences $(\al_k)\subset \cz$, $\eps_k=\pm 1$
and permutations $\pi$ of the positive integers. In this case, $X$
is called a symmetric space. The least constant $\la$ (over all
possible symmetric bases of $X$) is denoted by ${\rm sym}(X)$.
Again, we transfer this definition to the operator space setting:
$(x_k)$ is  completely symmetric if
  \[
  \|\sum_k \eps_k a_{\pi(k)} \ten x_k \|_p \le\la\,
  \|\sum_k  a_{k } \ten x_k \|_p
  \]
holds for finite sequences $(a_k)\subset S_p$, $\eps_k=\pm 1$ and
permutations $\pi$. If $X$ is a completely symmetric space, the
relevant constant is denoted by ${\rm sym}_{cb}(X)$. It is clear
that the four spaces in the previous corollary are completely
symmetric. Thus we deduce the following

\begin{cor}
 Let $\M$ and $p$ be as above. Then  every infinite dimensional subspace
of $L_p(\M)$ contains an infinite  completely symmetric basic
sequence.
 \end{cor}

It is not known whether the assertion above holds for $1\le p<2$.
This problem is open even for scalar coefficients. On the other
hand, we neither know whether the hyperfiniteness assumption can
be removed for $2<p<\infty$.  We refer to \cite{rx-JFA} and
\cite{ran-KP} for different versions of the Kadec-Pe\l czy\'{n}ski
alternative, which are most often at the Banach space level.

\smallskip

We now show that conversely all completely symmetric subspaces of
noncommutative $L_p$ are only those found in Corollary \ref{ssp}.
The next result is our starting point.

\begin{theorem}\label{cr}
 Let $\M$ be a von Neumann algebra, $2\le p<\8$  and $x_{ij}\in
L_p(\M)$. Then
 \[
   \big(\ez\, \big\|\sum_{i=1}^n  \eps_i x_{i\,\pi(i)} \big\|_p^p
 \big)^{1/p}
 \sim_{c_p}\max\big\{ \big(\frac 1n\, \sum_{i,j=1}^n
 \|x_{ij}\|_p^p \big)^{1/p}\,,\; \big\|\big(\frac 1n\,
\sum_{i,j=1}^n
 (x_{ij}^*x_{ij}+x_{ij}x_{ij}^*) \big)^{1/2}\big\|_p\big\} .
  \]
Here the expectation $\ez $ is taken over all choices of signs
$\eps_i=\pm1$ and all permutations $\pi$ on $\{1,...,n\}$.
\end{theorem}

\pf Again, we can assume that $\M$ is equipped with a normal
faithful state $\f$.  We consider $\Om=\{-1,1\}^n\times \Pi_n$,
where $\Pi_n$ is the set of all permutations on $\{1,...,n\}$. The
Haar measure on this group is the product  measure $\mu=\eps \ten
\nu$ of the normalized counting measures $\eps$ and $\nu$ on
$\{-1,1\}^n$, $\Pi_n$, respectively. The underlying von Neumann
algebra is then given by $(\N,\psi)=L_\8(\Om,2^{\Om},\mu)\ten
(\M,\f)$. In order to apply the noncommutative Burkholder
inequality we have to use the right filtration taken from
\cite{jmst}. For $k=1,...,n$ we consider the functions $f_k:
\Pi_n\to \rz$, $f_k(\pi)=\pi(k)$. The $\si$-algebra $\Si^2_k$ is
defined as the smallest $\si$-algebra on $\Pi_n$ making
$f_1,...,f_k$ measurable. By $\Si^1_k$ we denote the smallest
$\si$-algebra on $\{-1,1\}^n$ making $\eps_1,...,\eps_k$
measurable, where $\eps_1,...,\eps_n$ are the coordinate functions
on $\{-1,1\}^n$. Let $\Si_k$ be the product $\si$-algebra
$\Si^1_k\times \Si_k^2$. We then define the filtration $(\N_k)_k$
of $\psi$-invariant subalgebras  by
 \[ \N_k =L_\8(\Om,\Si_k,\mu) \ten \M \,. \]
Let $\ez_k$ be the conditional expectation associated to $\Si_k$.
Then $\E_k=\ez_k\ten {\rm id}$ is the state preserving conditional
expectation from $\N$ onto $\N_k$.

After these preliminaries, we consider
 \[ x =\sum_{i=1}^n \eps_i x_{i\,\pi(i)} \in L_p(\N) .\]
Let $d_k=dx_k$ be the martingale differences of $x$ with respect
to $(\N_k)$.  We note that $d_k=\eps_k x_{k\,\pi(k)}$. Therefore,
the noncommutative Burkholder inequality \cite{jx-burk} implies
 \[
 \|x\|_p\le c_p\max\big\{
 \big(\sum_{k=1}^n \|x_{k\,\pi(k)}\|_p^p \big)^{1/p}\,,\;
 \big\|\sum_{k=1}^n\E_{k-1}(x_{k\,\pi(k)}^* x_{k\,\pi(k)}
 + x_{k\,\pi(k)} x_{k\,\pi(k)}^*) \big\|_{p/2}^{1/2}\big\} .
 \]
Clearly, for every $k=1,...,n$, we have
 \be
  \|x_{k\,\pi(k)} \|_p^p
  &=&\sum_{j=1}^n \nu(\{\pi: \pi(k)=j\})\, \|x_{kj}\|_p^p\\
  &=&\sum_{j=1}^n \frac{(n-1)!}{n!}\,  \|x_{kj}\|_p^p
  =\frac{1}{n}\, \sum_{j=1}^n   \|x_{kj} \|_p^p  \, .
 \ee
Hence,
 \[ \big(\sum_{k=1}^n \|x_{k\,\pi(k)}\|_p^p \big)^{1/p}
 =\big(\frac{1}{n}\, \sum_{k,j=1}^n   \|x_{kj} \|_p^p\big)^{1/p}
 \, .\]
Let $\ez_{k}^2$ be the conditional expectation onto
$L_\8(\Pi_n,\Si_k^2,\mu)$. We observe that
 \[ \E_{k-1}(x_{k\,\pi(k)}^* x_{k\,\pi(k)}) =
  (\ez_{k-1}^2\ten{\rm id})( x_{k\,\pi(k)}^*
  x_{k\,\pi(k)}) \,.\]
The atoms in $\Si_{k-1}^2$ are indexed  by $(k-1)$-tuples
$(i_1,...,i_{k-1})$ of distinct integers in $\{1,...,n\}$. More
precisely,
 \[ A_{(i_1,...,i_{k-1})} =\big\{ \pi\;:\;  \pi(1)=i_1, ...,
 \pi(k-1)= i_{k-1}\big\} .\]
Clearly, the cardinality of $A_{(i_1,...,i_{k-1})}$ is that of
$\Pi_{n-(k-1)}$, i.e., $(n-k+1)!$. Therefore, letting
$\al_k=(n-k+1)!/n!$, we get
 \be
 (\ez_{k-1}^2\ten{\rm id})( x_{k\,\pi(k)}^* x_{k\,\pi(k)})=
 \sum_{(i_1,...,i_{k-1})} \un_{A_{(i_1,...,i_{k-1})}}\;  \al_k^{-1}
 \int_{A_{(i_1,...,i_{k-1})}} x_{k\,\pi(k)}^*x_{k\,\pi(k)}
 \, d\nu(\pi)\,.
 \ee
For fixed $(i_1,...,i_{k-1})$, letting $B=\{i_1,...,i_{k-1}\}$, we
have
 \[
 \al_{k}^{-1}  \int_{A_{(i_1,...,i_{k-1})}}
 x_{k\,\pi(k)}^*x_{k\,\pi(k)} \, d\nu(\pi)
 =\frac{1}{n-k+1}\, \sum_{j\notin B} x_{kj}^*x_{kj} \,.
 \]
Hence for all $k\le n/2$ we deduce
 \[ (\ez_{k-1}^2\ten{\rm id})( x_{k\,\pi(k)}^* x_{k\,\pi(k)})
 \le  \frac{2}{n}\,\sum_{j=1}^n x_{kj}^*x_{kj} \,.\]
Let us assume temporarily that $x_{kj}=0$ for $k> n/2$. Then
combining the previous estimates, we obtain
 \[ \sum_{k=1}^n
  \E_{k-1}(x_{k\,\pi(k)}^* x_{k\,\pi(k)}) \le
  \frac{2}{n}\, \sum_{k,j=1}^n x_{kj}^*x_{kj} \]
for all permutations $\pi$. The same argument applies to
$x_{k\,\pi(k)}x_{k\,\pi(k)}^*$ too. Therefore, we get  the upper
estimate under the additional assumption that $x_{kj}=0$ for $k>
n/2$. The general case then follows from triangle inequality.

For the lower estimate we use the Jensen inequality  and the
orthogonality of the Rademacher variables (noting that $p/2\ge
1$):
 \be
 \|x\|_p^2
 &=&\ez\, \|x^*x \|_{p/2}\ge \|\ez(x^*x)\|_{p/2}\\
 &=&\big\|\sum_{k=1}^n \int_{\Pi_n} x_{k\,\pi(k)}^*x_{k\,\pi(k)}
  \big\|_{p/2}
  = \big\|\frac{1}{n}\,\sum_{k,j=1}^n x_{kj}^*x_{kj}\big\|_{p/2} \,.
 \ee
The same calculation  involving $xx^*$ yields the other square
function estimate. Since $L_p(\N)$ has cotype $p$, we easily find
the missing estimate on the diagonal term. \qd


\begin{cor}\label{symn}
 Let  $\A$ and $\M$  be von Neumann algebras and   $2\le p<\infty$.
Let  $(x_k)_{1\le k\le n}\subset L_p(\M)$  and $\la>0$ such that
 \[
 \big\|\sum_{k=1}^n \eps_k a_{\pi(k)} \ten  x_k\big\|_p \le\la\,
 \big\|\sum_{k=1}^n a_{k} \ten  x_k\big\|_p \]
holds for all $\eps_k=\pm 1$, all  permutations $\pi$ on
$\{1,...,n\}$ and coefficients $a_k\in L_p(\A)$. Then there are
constants $\al,\beta$ and $\gamma$, depending only on $(x_k)$,
such that
 \[
 \big\|\sum_{k=1}^n a_{k} \ten  x_k\big\|_p\sim_{\la^2c_p}
 \max\big\{\al\big(\sum_{k=1}^n \|a_k\|_p^p\big)^{1/p}\,,\;
 \beta \big\|(\sum_{k=1}^n a_k^*a_k)^{1/2}\big\|_p\,,\;
 \gamma \big\|(\sum_{k=1}^n a_ka_k^*)^{1/2}\big\|_p\big\}\,.
 \]
holds for all $a_k\in L_p(\A)$.
\end{cor}

\pf This is an easy consequence of Theorem \ref{cr}. Indeed, we
have
 \[ \frac{1}{\la}\,\big\|\sum_{k=1}^n a_{k}\ten  x_k\big\|_p
 \le\big(\ez \big\|\sum_{k=1}^n \eps_k a_{\pi(k)} \ten
 x_k\big\|_p^p\big)^{1/p}
 \le \la \big\|\sum_{k=1}^n a_{k} \ten  x_k\big\|_p \,.\]
Then we deduce the assertion with
 \[
  \al= \big(\frac1n\,\sum_{k=1}^n \|x_k\|_p^p\big)^{1/p}\,,\quad
 \beta= \big\|(\frac1n\,\sum_{k=1}^n x_k^*x_k)^{1/2}\big\|_p\,,\quad
 \gamma=\big\|(\frac1n\,\sum_{k=1}^n x_kx_k^*)^{1/2}\big\|_p \,.
 \] \qd

\smallskip

Let us introduce a more notation. For a Banach (or operator) space
$X$ and positive real $\al$, $\al X$ denotes $X$ but equipped with
the norm $\al\|\,\|$. For convenience, set $\al X=\{0\}$ if
$\al=0$. Recall that the Banach-Mazur distance between two Banach
spaces $X$ and $Y$ is
 \[ d(X,\,Y) = \inf\big\{\|T\|\,\|T^{-1}\| :\;\;
 T:E\to F\; \mbox{isomorphism}\big\}.\]
Similarly, we define the operator space analogue $d_{cb}(X, Y)$ by
replacing the norm of an isomorphism by the cb-norm of a complete
isomorphism.

\begin{cor}\label{fd} Let $2\le p<\8$ and $X$ be an $n$-dimensional
subspace of $L_p(\M)$. Then
 \begin{enumerate}[\rm(i)]
 \item there exist nonnegative $\al$ and
$\beta$ such that
 \[ d(X,\; \al\ell_p^n\cap \beta\ell_2^n)
 \le c_p \, {\rm sym}(X)^2 \, ;\]
 \item there exist nonnegative $\al$,
$\beta$ and $\gamma$ such that
 \[ d_{cb}(X,\; \al\ell_p^n\cap \beta C_p^n
 \cap \gamma R_p^n)\le c_p \, {\rm sym}_{cb}(X)^2 \, .\]
 \end{enumerate}
\end{cor}

\pf Let $(x_1,...,x_n)$  be a (completely) symmetric basis of $X$
with constant $\la\le 2\,{\rm sym}(X)$ (or $\la\le 2\,{\rm
sym}_{cb}(X)$). Such a basis exists for $\dim X<\8$. It then
remains to apply the previous corollary with $\A=\cz$ for (i) and
$\A=B(\ell_2)$ for (ii). \qd

\begin{cor}\label{symmfin}
Let $2\le p<\infty$ and $X\subset L_p(\M)$ be an infinite
dimensional subspace.
 \begin{enumerate}[\rm(i)]
  \item If $X$ is symmetric, then $X$ is isomorphic to $\ell_p$ or
$\ell_2$.
 \item If $X$ is completely symmetric,
then $X$ is completely isomorphic to $\ell_p$, $C_p$, $R_p$ or
$C_p\cap R_p$.
 \end{enumerate}
\end{cor}

\pf We prove only (ii). The proof of (i) is simpler, just by
replacing vector coefficients by scalar ones. Let $(x_k)$ be a
completely symmetric basis of $X$ with constant $\la$. For every
$\nen$,  set
  \[
  \al_n= \big(\frac1n\,\sum_{k=1}^n \|x_k\|_p^p\big)^{1/p}\,,\quad
 \beta_n= \big\|(\frac1n\,\sum_{k=1}^n x_k^*x_k)^{1/2}\big\|_p\,,\quad
 \gamma_n=\big\|(\frac1n\,\sum_{k=1}^n x_kx_k^*)^{1/2}\big\|_p \, .
 \]
Note that $\al_n$, $\beta_n$ and $\gamma_n$ are less than or equal
to $\sup_k \|x_k\|_p$. Passing to subsequences if necessary, we
may assume that the three sequences $(\al_n)$, $(\beta_n)$ and
$(\gamma_n)$ converge respectively to $\al$, $\beta$ and $\gamma$.
Thus by Corollary \ref{symn}, for any finite  sequence
$(a_k)\subset S_p$ we have
 \[
 \big\|\sum_{k} a_{k} \ten  x_k\big\|_p\sim_{\la^2c_p}
 \max\big\{\al\big(\sum_{k}\|a_k\|_p^p\big)^{1/p}\,,\;
 \beta \big\|(\sum_{k} a_k^*a_k)^{1/2}\big\|_p\,,\;
 \gamma \big\|(\sum_{k} a_ka_k^*)^{1/2}\big\|_p\big\}\,.
 \]
Then using \eqref{p-2}, we deduce that $X$ is completely
isomorphic to $\ell_p$ if $\beta=\gamma=0$, to $C_p$ if $\beta>0$
and $\gamma=0$, to $R_p$ if $\beta=0$ and $\gamma>0$, and finally
to $R_p\cap C_p$ if $\beta>0 $ and $\gamma>0$.\qd

\begin{rem}{\rm
 It will be shown in \cite{juros} that every subsymmetric
basic sequence in $L_p(\M)$ is symmetric. A sequence $(e_k)$ is
called subsymmetric if
 \[ \|\sum_k \eps_k a_k\ten e_{j_k} \|_p \sim_c \|\sum_k
 a_k\ten e_k\| \]
holds for all increasing sequences $(j_k)$ of integers. Therefore,
Corollary \ref{symmfin} yields a characterization of subspaces of
noncommutative $L_p$ with a subsymmetric basis.}
\end{rem}

If $\M$ is finite, we can eliminate the two spaces $C_p$ and $R_p$
in Corollary \ref{symmfin} (ii).

\begin{cor} Let $2<p<\infty$ and $\M$ be a finite von Neumann
algebra. Then  $C_p$ and $R_p$ do not completely embed into
$L_p(\M)$.
\end{cor}

\pf We assume that $\f$ is a normal faithful tracial state on
$\M$. Suppose that $C_p$ completely embeds into $L_p(\M)$. Namely,
there exists an infinite sequence $(x_k)\subset L_p(\M)$ such that
 \[\big\|\sum_k a_k\ten x_k\big\|_p \sim
 \big\|\big(\sum_k a_k^*a_k\big)^{1/2}\big\|_p\]
holds for all $(a_k)\subset S_p$. In particular, if
$\al=(\al_{ik})\in S_p$, then
 \[ \|\al\|_{S_p} \sim
 \big\|\sum_{i,k} \al_{ik}e_{1\,i}\ten x_k
 \big\|_{L_p(B(\ell_2)\bar{\ten}\M)}\,.  \]
Note  that for
 \[x=\sum_{i,k} \al_{ik} e_{1\,i}\ten x_k\quad
 \mbox{and}\quad x_i=\sum_k \al_{ik}x_k\,,
 \]
the H\"{o}lder inequality implies
 \be
  \|x\|_{L_2(B(\ell_2)\bar{\ten}\M)}^2
  &=& \|xx^*\|_{L_1(B(\ell_2)\bar{\ten}\M)}
  = \|\sum_i x_ix_i^*\|_{L_1(\M)}\\
  &\le& \|\sum_i x_ix_i^*\|_{L_{p/2}(\M)}
  = \|x\|_{L_p(B(\ell_2)\bar{\ten}\M)}^2 \,.
  \ee
Here the $L_p(\M)$ are defined in terms of the trace $\f$. This
tells us that on the subspace $Y={\rm span}\{ e_{1\,i}\ten x_k\}$
the norms in $L_p\cap L_2$ and $L_p$ coincide. Thus we have found
an embedding of $S_p$ into $L_p(B(\ell_2)\bar{\ten}\M)\cap
L_2(B(\ell_2)\bar{\ten}\M)$. According to \cite{ju-poisson} the
latter space embeds into $L_p(\R)$ for a finite von Neumann
algebra $\R$. Thus we obtain an embedding of $S_p$ into $L_p(\R)$.
This is, however, absurd in view of the results in
\cite{suk-noniso}. \qd


\section{Bisymmetric and unitary invariant subspaces of
$L_p$}
 \label{unitary ideal}


We extend in this section the results in the previous one to the
case of double indices. Namely, we will determine the bisymmetric
and unitary invariant subspaces of noncommutative $L_p$-spaces for
$2<p<\8$. In particular,  we will characterize those unitary
ideals which can embed into a noncommutative $L_p$. For notational
convenience, given a finite matrix $x=(x_{ij})$ with entries in
$L_p(\M)$ we introduce
 \begin{align*}
 \gamma_0(x)
 &=\big(\sum_{i,j}\|x_{ij}\|_p^p\big)^{1/p}\,,\\
 \gamma_1(x)
 &=\big(\sum_{i}\big\|\big(\sum_{j}x_{ij}^*x_{ij}\big)^{1/2}
 \big\|_p^p\big)^{1/p}\,,\quad
 \gamma_2(x)
 =\big(\sum_{i}\big\|\big(\sum_{j}x_{ij}x_{ij}^*\big)^{1/2}
 \big\|_p^p\big)^{1/p}\,,\\
 \gamma_3(x)
 &=\big(\sum_{j}\big\|\big(\sum_{i}x_{ij}^*x_{ij}\big)^{1/2}
 \big\|_p^p\big)^{1/p}\,,\quad
 \gamma_4(x)
 =\big(\sum_{j}\big\|\big(\sum_{i}x_{ij}x_{ij}^*\big)^{1/2}
 \big\|_p^p\big)^{1/p}\,,\\
 \gamma_5(x)
 &=\big(\big\|\big(\sum_{i,j}x_{ij}^*x_{ij}\big)^{1/2}
 \big\|_p\,,\quad\hskip 1.4cm
 \gamma_6(x)
 =\big(\big\|\big(\sum_{i,j}x_{ij}x_{ij}^*\big)^{1/2}
 \big\|_p\,,\\
 \gamma_7(x)
 &=\big\|\sum_{i,j}e_{ij}\ten x_{ij}
 \big\|_{L_p(B(\ell_2)\bar\ten\M))}\,,\quad\hskip 0.1cm
 \gamma_8(x)
 =\big\|\sum_{i,j}e_{ji}\ten x_{ij}
 \big\|_{L_p(B(\ell_2)\bar\ten\M))}\,.
 \end{align*}

\begin{theorem} \label{bs}
 Let $2\le p<\8$ and  $\A$ and $\M$ be von Neumann algebras.
Let  $a=(a_{ij})$  and $x=(x_{ij})$ be two $n\times n$ matrices
with entries in  $L_p(\A)$ and $L_p(\M)$, respectively.  Then
 \begin{align*}
 &\big(\ez\,\big\|\sum_{i,j=1}^n\eps_i\eps_j'\,a_{ij}\ten
 x_{\pi(i)\,\pi'(j)}\big\|_p\big)^{1/p}
  \sim_{c_p} \\
 &\frac1{n^{2/p}}\,\gamma_0(a)\gamma_0(x)
  + \frac1{n^{1/p +1/2}}\,\sum_{k=1}^4\gamma_k(a)\gamma_k(x)
 +\frac1n\,\sum_{k=5}^8\gamma_k(a)\gamma_k(x)\,.
 \end{align*}
Here the expectation $\ez$ is taken over independent copies
$\eps_i$, $\eps'_i$ of Rademacher variables  and independent
copies $\pi$ and $\pi'$  of permutations on $\{1,...,n\}$.
\end{theorem}

\pf This is an iteration of Theorem  \ref{cr}. By that theorem,
we get
  \begin{align*}
  & \big(\ez\,\big\|\sum_{i,j=1}^n\eps_i\eps_j'\,a_{ij}\ten
  x_{\pi(i)\,\pi'(j)}\big\|_p\big)^{1/p}
   \sim \\
  &\hskip 1.5cm  \frac1{n^{1/p}}\,\big(\ez_{\eps', \pi'}\,\sum_{i,k}
   \big\|\sum_{j}\eps_j'\,a_{ij}\ten
  x_{k\,\pi'(j)}\big\|_p^p\big)^{1/p}\\
  &\hskip 1cm +  \frac1{n^{1/2}}\, \big(\ez_{\eps', \pi'}\,\big\|\sum_{i,k,j}
   \eps_j'\,e_{(i,k),(1,1)}\ten a_{ij}\ten
  x_{k\,\pi'(j)}\big\|_p^p\big)^{1/p}\\
  &\hskip 1cm +  \frac1{n^{1/2}}\,\big(\ez_{\eps', \pi'}\,\big\|\sum_{i,k,j}
   \eps_j'\,e_{(1,1),(i,k)}\ten a_{ij}\ten
  x_{k\,\pi'(j)}\big\|_p^p\big)^{1/p}\,
  \mathop{=}^{\rm def}\, I+II+III\,.
  \end{align*}
Here we use $e_{(i,k),(j,l)}$ to denote the matrix units of
$B(\ell_2(\nz^2))$, so $(i,k)$ and $(j,l)$ index rows and columns,
respectively.  We apply Theorem \ref{cr} for a second time to the
first term on the right hand side and find
 \be
  I
  &\sim& \frac1{n^{2/p}}\, \gamma_0(a)\gamma_0(x)
  +\frac1{n^{1/p +1/2}}\,\big(\sum_{i,k} \big\|\sum_{j,l}
  e_{(j,l),(1,1)} \ten a_{ij}\ten x_{kl}  \big\|_p^p \big)^{1/p}\\
  &&+ \;\frac1{n^{1/p +1/2}}\, \big(\sum_{i,k} \big\|\sum_{j,l}
  e_{(1,1),(j,l)} \ten a_{ij}\ten x_{kl}  \big\|_p^p \big)^{1/p}
    \,.
 \ee
Identifying $e_{(j,l),(1,1)}$ with $e_{j1}\ten e_{l1}$ (up to a
conjugation by a unitary), we have
 \be
  \big\|\sum_{j,l}
  e_{(j,l),(1,1)} \ten a_{ij}\ten x_{kl}  \big\|_p
   &=& \big\|\big(\sum_{j} e_{j1} \ten a_{ij}\big)\ten
   \big(\sum_{l} e_{l1} \ten x_{kl}\big) \big\|_p\\
 &=& \big\|\big(\sum_{j}a_{ij}^* a_{ij}\big)^{1/2}\big\|_p\,
   \big\|\big(\sum_{l}x_{kl}^* x_{kl}\big)^{1/2}\big\|_p \, .
 \ee
We deal with similarly the other term containing $e_{(1,1),(j,l)}$
and then deduce that
 \[ I\sim \frac1{n^{2/p}}\,\gamma_0(a)\gamma_0(x)
 + \frac1{n^{1/p +1/2}}\, \gamma_1(a)\gamma_1(x)
 + \frac1{n^{1/p +1/2}}\, \gamma_2(a)\gamma_2(x)
 \,.\]
Similar arguments apply to $II$ and $III$ too. $II$ is again
equivalent to a sum of three terms. Let us consider, for instance,
the second one on column norm, which is
  \begin{align*}
  &\frac1n\,\big\|\sum_{j,l}\sum_{i,k}e_{(j,l),(1,1)}\ten
  e_{(i,k),(1,1)} \ten a_{ij}\ten x_{kl}\big\|_p\\
  &=\frac1n\,\big\|\sum_{i,j}e_{(i,j),(1,1)}\ten a_{ij}\big\|_p\,
  \big\|\sum_{k,l}e_{(k,l),(1,1)}\ten x_{kl}\big\|_p
  =\frac1n\,\gamma_5(a)\gamma_5(x)\,.
  \end{align*}
Then we see that
 \[II\sim \frac1{n^{1/p +1/2}}\,\gamma_3(a)\gamma_3(x)
 +\frac1n\,\gamma_5(a)\gamma_5(x) +
 \frac1n\,\gamma_7(a)\gamma_7(x)\,.\]
Finally, $III$ yields the three missing terms. \qd

\smallskip

Permutations and $(\eps_1,..., \eps_n)$ induce permutation and
diagonal matrices, which are, of course, unitary. If the
expectation in Theorem \ref{bs} is taken over all unitary
matrices, we get a much simpler equivalence.

\begin{theorem}\label{unitary}
Under the assumption of Theorem \ref{bs}, we have
 \begin{align*}
 \big(\ez\,\big\|\sum_{i,j, k, l=1}^n
 u_{ik}v_{lj}\,a_{ij}\ten x_{kl}\big\|_p\big)^{1/p}
  \,\sim_{c_p}\,
 \frac1n\,\sum_{k=5}^8\gamma_k(a)\gamma_k(x)\,.
 \end{align*}
Here the expectation $\ez$ is the integration in $(u_{ik})$ and
$(v_{lj})$ on $U(n)\times U(n)$, where $U(n)$ is the $n\times n$
unitary group equipped with Haar measure.
\end{theorem}

\pf The proof is similar to that of Theorem \ref{bs}. Instead of
the noncommutative Burkholder inequality via Theorem \ref{cr}, we
now use the noncommutative Khintchine inequality with help of the
classical fact that $(u_{ik})$ can be replaced by a Gaussian
matrix $n^{-1/2}\,(g_{ij})$, where the $g_{ij}$ are independent
Gaussian variables of mean-zero and variance $1$ (see
\cite{mar-pis}). Thus
 \[\big(\ez\,\big\|\sum_{i,j, k, l=1}^n
 u_{ik}v_{lj}\,a_{ij}\ten x_{kl}\big\|_p\big)^{1/p}
  \,\sim_{c}\,
  \frac1n\,\big(\ez\,\big\|\sum_{i,j, k, l=1}^n
 g_{ik}\,g'_{lj}\,a_{ij}\ten x_{kl}\big\|_p\big)^{1/p}
 \,.\]
It then remains to repeat the arguments in the proof of Theorem
\ref{bs} by using
 \[\big(\ez\,\big\|\sum_{i, k}
 g_{ik}\, y_{ik}\big\|_p\big)^{1/p}\,\sim_{c\sqrt p}\,
 \big\|\big(\sum_{i, k}y_{ik}^*y_{ik}\big)^{1/2}\big\|_p
 +\big\|\big(\sum_{i, k}y_{ik}y_{ik}^* \big)^{1/2}\big\|_p\]
for any $y_{ik}$ in a noncommutative $L_p$ (see \cite{pis-ast}).
\qd

\smallskip

We say that $(x_{ij})$ is a bisymmetric basis of a subspace
$X\subset L_p(\M)$ if every entry $x_{ij}$ is nonzero, the linear
span of the $x_{ij}$ is dense in $X$ and there exists a constant
$\la$ such that
 \[\big\|\sum_{i,j}\eps_i\eps'_j a_{\pi(i)\,\pi'(j)} x_{ij}\big\|_p
 \le\la\,\big\|\sum_{i,j}a_{ij}x_{ij}\big\|_p \]
holds for all finite scalar matrices $(a_{ij})$, all $\eps_i=\pm
1$, $\eps'_j=\pm 1$ and all permutations $\pi$ and $\pi'$. It is
easy to check that $(x_{ij})$ is indeed a basis of $X$ according
to an appropriate order, for instance, the one defined as follows.
Let $e_1=x_{11}$ and assume defined $e_1, ..., e_{n^2}$. Then we
set $e_{n^2+j}=x_{n+1,j}$ for $j=1,..., n+1$, and
$e_{n^2+n+1+i}=x_{n+1-i,n+1}$ for $i=1,..., n$. Similarly, we
define completely bisymmetric bases by replacing scalar
coefficients by matrices.

\smallskip

Recall that $\ell_p(\ell_2)$ denotes the space of all scalar
matrices $a=(a_{ij})$ such that
 \[\big(\sum_{i}\big(\sum_{j}|a_{ij}|^2\big)^{p/2}
 \big)^{1/p}<\8\]
and is equipped with the natural norm.  $\ell_p(\ell_2)^{t}$ is
the space of all $a$ such that $a^t\in \ell_p(\ell_2)$, where
$a^t$ denotes the transpose of $a$.  In the operator space
setting, these spaces yield four different spaces  $\ell_p(C_p),
\ell_p(R_p), \ell_p(C_p)^{t}$ and $\ell_p(R_p)^{t}$, corresponding
respectively to the norms $\gamma_k$ for $1\le k\le 4$ introduced
at the beginning of the present section. Accordingly, the last
four norms there give four other operator spaces $ C_p(\nz^2),
R_p(\nz^2), S_p$ and $S_p^t$. The following is the matrix analogue
of Corollary \ref{symmfin}. The proof is almost identical to that
of Corollary \ref{symmfin} but now via Theorem \ref{bs}.

\begin{cor}\label{bisym}
Let $2\le p<\infty$ and $X\subset L_p(\M)$ be a subspace.
 \begin{enumerate}[\rm(i)]
  \item If $X$ has a bisymmetric basis given by an infinite
matrix, then $X$ is isomorphic to one of the following spaces
 \[ \ell_p(\nz^2),\;  \ell_p(\ell_2),\; \ell_p(\ell_2)^{t},\;
 \ell_p(\ell_2)\cap\ell_p(\ell_2)^{t},\; S_p,\;\ell_2(\nz^2)\,.
 \]
 \item If $X$ has a completely bisymmetric basis given by an infinite
matrix, then $X$ is completely isomorphic to one of the following
spaces
 \[ \ell_p(\nz^2),\;  \ell_p(C_p),\; \ell_p(R_p),\;
  \ell_p(C_p)^{t},\;\ell_p(R_p)^{t},\; S_p,\; {S_p}^t,\;
   C_p(\nz^2),\; R_p(\nz^2)
 \]
or one possible intersection of them.
 \end{enumerate}
\end{cor}

\begin{rem} {\rm
 The relations between the $9$ building blocks in (ii) above
are shown by the following diagram
 
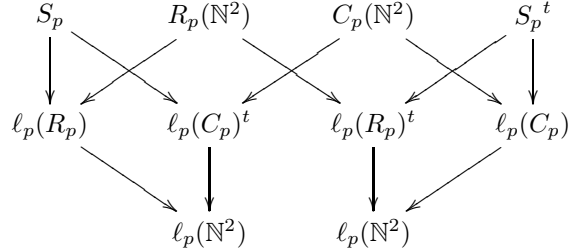
\[\xymatrix{
 & S_p \ar[d] \ar[rd] & R_p(\nz^2) \ar[ld] \ar[rd]
 & C_p(\nz^2)\ar[ld] \ar[rd] & {S_p}^{t}\ar[ld] \ar[d] \\
 & \ell_p(R_p) \ar[rd] & \ell_p(C_p)^t\ar[d]
 & \ell_p(R_p)^{t}\ar[d] &  \ell_p(C_p) \ar[ld] \\
 &&\ell_p(\nz^2)&\ell_p(\nz^2) &
 }\]
The arrows indicate complete contractions, e.g.,  $S_p\subset
\ell_p(R_p)\cap \ell_p(C_p)^t$ and $C_p(\nz^2) \subset
\ell_p(C_p)\cap \ell_p(C_p)^{t}$. Thus not all intersections of
these spaces are nontrivial for some of them simplify. However,
the four spaces on each of the first two levels do not give any
nontrivial intersection, so yield 16 pairwise distinct spaces.
 }\end{rem}

\begin{rem} {\rm
 It is easy to see that all spaces appearing the preceding
corollary (completely) embed really into a noncommutative $L_p$.
Note that an interesting embedding of
$\ell_p(\ell_2)\cap\ell_p(\ell_2)^{t}$ (or $\ell_p(R_p)\cap
  \ell_p(C_p)^{t}$ in the operator space case) is given by the
  noncommutative Khintchine inequality in Remark \ref{random Sp}.
 }\end{rem}

A bysymmetric basis $(x_{ij})$  of $X$ is called (completely)
unitary invariant if
 \[\big\|\sum_{ij}u_{ij}a_{ij}x_{ij}\big\|_p\le\la\,
 \big\|\sum_{ij}a_{ij}x_{ij}\big\|_p\]
holds for all unitaries $(u_{ij})$ and all finite matrices
$(a_{ij})$ in $\cz$ (in $S_p$). Recall that if $E$ is a symmetric
sequence space, the associated unitary ideal $S_E$ is defined to
be the closure of finite matrices with respect to the norm
 \[ \|a\|_{S_E} = \|(s_k(a))_k \|_E \, ,\]
where $(s_k(a))_k$ is the sequence of the singular numbers of $a$.
It is well known that the matrix units of $B(\ell_2)$ form a
unitary invariant basis of $S_E$.

\begin{cor}\label{unitary inv}
 Let $2\le p<\infty$ and $X\subset L_p(\M)$ be a subspace.
 \begin{enumerate}[\rm(i)]
  \item If $X$ has a unitary invariant basis,
 then $X$ is isomorphic to $S_p$ or $S_2$. Consequently, a unitary
 ideal $S_E$ embeds in $L_p(\M)$ iff $E=\ell_p$ or $E=\ell_2$.
  \item If $X$ has a completely unitary invariant basis,
then $X$ is completely isomorphic to one of the $16$ spaces:
 $S_p,\; {S_p}^t,\; C_p(\nz^2),\; R_p(\nz^2)$
and their intersections.
 \end{enumerate}
 \end{cor}

\pf This is an immediate consequence of the preceding corollary
since all spaces there but those in the present corollary are not
unitary invariant. Alternately, we can also follow the proof of
Corollary \ref{bisym} by using Theorem \ref{unitary}.\qd


\begin{thebibliography}{JNRX}

\bibitem[BL]{bl}
J.~Bergh, J.~L{\"o}fstr{\"o}m.
\newblock {\em Interpolation spaces.}
\newblock Springer-Verlag, Berlin, 1976.

\bibitem[BKS]{boks}
M.~Bo{\.z}ejko, B.~K{\"u}mmerer,  R.~Speicher.
\newblock {$q$}-{G}aussian processes: non-commutative and classical aspects.
\newblock {\em Comm. Math. Phys.}, 185:129--154, 1997.

\bibitem[BS1]{bos-example}
M.~Bo{\.z}ejko, R.~Speicher.
\newblock An example of a generalized {B}rownian motion.
\newblock {\em Comm. Math. Phys.}, 137:519--531, 1991.

\bibitem[BS2]{bos-cox}
M.~Bo{\.z}ejko, R.~Speicher.
\newblock Completely positive maps on {C}oxeter groups, deformed commutation
  relations, and operator spaces.
\newblock {\em Math. Ann.}, 300:97--120, 1994.


\bibitem[Bu]{burk-dis}
D.~L.~Burkholder.
\newblock Distribution function inequalities for martingales.
\newblock {\em Ann. Probability}, 1:19--42, 1973.

\bibitem[BuG]{burk-gundy}
D.~L.~Burkholder, R.~F.~Gundy.
\newblock Extrapolation and interpolation of quasi-linear operators on
  martingales.
\newblock {\em Acta Math.}, 124:249--304, 1970.


\bibitem[C]{con-classi}
A.~Connes.
\newblock Classification of injective factors. {C}ases {$II\sb{1},$}
  {$II\sb{\infty },$} {$III\sb{\lambda },$} {$\lambda \not=1$}.
\newblock {\em Ann. Math.}, 104:73--115, 1976.

\bibitem[ER]{er-book}
Ed.~Effros, Z-J.~Ruan.
\newblock {\em Operator spaces}, volume~23 of {\em London Mathematical Society
  Monographs. New Series}.
\newblock The Clarendon Press Oxford University Press, New York, 2000.

\bibitem[GGMS]{ggms}
  N.~Ghoussoub, G.~Godefroy, B.~Maurey,  W.~Schachermayer.
\newblock Some topological and geometrical structures in {B}anach spaces.
\newblock {\em Mem. Amer. Math. Soc.}, 70(378):iv+116, 1987.

\bibitem[H1]{haag-Lp}
U.~Haagerup.
\newblock {$L\sp{p}$}-spaces associated with an arbitrary von {N}eumann
  algebra.
\newblock In {\em Alg\`ebres d'op\'erateurs et leurs applications en physique
  math\'ematique (Proc. Colloq., Marseille, 1977)}, volume 274 of {\em Colloq.
  Internat. CNRS}, pages 175--184. CNRS, Paris, 1979.

\bibitem[H2]{haag-bicent}
U.~Haagerup.
\newblock Connes' bicentralizer problem and uniqueness of the injective factor
  of type {${\rm III}\sb 1$}.
\newblock {\em Acta Math.}, 158:95--148, 1987.


\bibitem[HW]{haag-win}
U.~Haagerup, C.~Winslow.
\newblock The {E}ffros-{M}ar\'echal topology in the space of von {N}eumann
  algebras. {II}.
\newblock {\em J. Funct. Anal.}, 171:401--431, 2000.

\bibitem[Hi]{hiai}
F.~Hiai.
\newblock {$q$}-deformed {A}raki-{W}oods algebras.
\newblock In {\em Operator algebras and mathematical physics (Constan\c ta,
  2001)}, pages 169--202. Theta, Bucharest, 2003.

\bibitem[JMST]{jmst}
W.~B.~Johnson, B.~Maurey, G.~Schechtman, L.~Tzafriri.
\newblock Symmetric structures in {B}anach spaces.
\newblock {\em Mem. Amer. Math. Soc.}, 19(217):v+298, 1979.

\bibitem[JSZ]{jsz}
 W.~B.~Johnson, G.~Schechtman, J.~Zinn.
\newblock Best constants in moment inequalities for linear combinations of
  independent and exchangeable random variables.
\newblock {\em Ann. Probab.}, 13:234--253, 1985.

\bibitem[J1]{ju-doob}
M.~Junge.
\newblock Doob's inequality for non-commutative martingales.
\newblock {\em J. Reine Angew. Math.}, 549:149--190, 2002.

\bibitem[J2]{ju-OH}
M.~Junge.
\newblock Embedding of the operator space {$OH$} and the logarithmic `little
  {G}rothendieck inequality'.
\newblock {\em Invent. Math.}, 161:225--286, 2005.

\bibitem[J3]{ju-araki}
M.~Junge.
\newblock Operator spaces and Araki-Woods factors --A quantum probabilistic
approach--.
\newblock {\em Int. Math. Res. Pap.}, ID 76978, pp. 87, 2006.

\bibitem[J4]{ju-poisson}
M.~Junge.
\newblock A noncommuttaive Poisson process.
\newblock In preparation.


\bibitem[JNRX]{jnrx-OLp}
M.~Junge, N.~J.~Nielsen, Z-J.~Ruan,  Q.~Xu.
\newblock {$\mathcal{COL}_p$} spaces---the local structure of non-commutative
  {$L\sb p$} spaces.
\newblock {\em Adv. Math.}, 187:257--319, 2004.

\bibitem[JP]{jupa-amalg}
M.~Junge, J.~Parcet.
\newblock Theory of amalgamated $L_p$-spaces in noncommutative probability.
\newblock Preprint 2006.

\bibitem[JPX]{jpx}
M.~Junge, J.~Parcet, Q.~Xu.
\newblock Rosenthal type inequalities for free chaos.
\newblock {\em Ann. Proba.}, to appear.

\bibitem[JR]{juros}
M.~Junge, H.~Rosenthal.
\newblock Noncommutative $L_p$-spaces are asymptotically stable.
\newblock In preparation.

\bibitem[JRX]{jrx-OLp2}
M.~Junge, Z-J.~Ruan, Q.~Xu.
\newblock Rigid {$\mathcal{OL}_p$} structures of non-commutative {$L\sb
  p$}-spaces associated with hyperfinite von {N}eumann algebras.
\newblock {\em Math. Scand.}, 96:63--95, 2005.

\bibitem[JX1]{jx-burk}
M.~Junge, Q.~Xu.
\newblock Noncommutative {B}urkholder/{R}osenthal inequalities.
\newblock {\em Ann. Probab.}, 31:948--995, 2003.

\bibitem[JX2]{jx-const}
M.~Junge, Q.~Xu.
\newblock On the best constants in some non-commutative martingale
  inequalities.
\newblock {\em Bull. London Math. Soc.}, 37:243--253, 2005.

\bibitem[JX3]{jx-erg}
M.~Junge, Q.~Xu.
\newblock Noncommutative maximal ergodic theorems.
\newblock {\em J. Amer. Math. Soc.}, 20:385-439, 2007.

\bibitem[KR]{kar-II}
R.~V.~Kadison, J.~R.~Ringrose.
\newblock {\em Fundamentals of the theory of operator algebras. {V}ol. {II}},
  volume~16 of {\em Graduate Studies in Mathematics}.
\newblock American Mathematical Society, Providence, RI, 1997.

\bibitem[LT]{LT-I}
J.~Lindenstrauss, L.~Tzafriri.
\newblock {\em Classical {B}anach spaces. {I}}.
\newblock Springer-Verlag, Berlin, 1977.

\bibitem[LP]{lust-khin}
F.~Lust-Piquard.
\newblock In\'egalit\'es de {K}hintchine dans {$C\sb p\;(1<p<\infty)$}.
\newblock {\em C. R. Acad. Sci. Paris}, 303:289--292, 1986.

\bibitem[LPP]{LPP}
F.~Lust-Piquard, G.~Pisier.
\newblock Noncommutative {K}hintchine and {P}aley inequalities.
\newblock {\em Ark. Mat.}, 29:241--260, 1991.


\bibitem[MP]{mar-pis}
M.~B.~Marcus, G.~Pisier.
\newblock {\em Random Fourier series with applications to harmonic
analysis.}
\newblock {\em Annals of Mathematics Studies}, 101.
Princeton University Press, Princeton, N.J., 1981.

\bibitem[PaR]{parcet-ran-gundy}
J.~Parcet, N.~Randrianantoanina.
\newblock Gundy's decomposition for non-commutative martingales and
  applications.
\newblock {\em Proc. London Math. Soc.}, 93:227--252, 2006.

\bibitem[P1]{pis-ast}
G.~Pisier.
\newblock Non-commutative vector valued {$L\sb p$}-spaces and completely
  {$p$}-summing maps.
\newblock {\em Ast\'erisque}, (247):vi+131, 1998.

\bibitem[P2]{pis-intro}
G.~Pisier.
\newblock {\em Introduction to operator space theory}, volume 294 of {\em
  London Mathematical Society Lecture Note Series}.
\newblock Cambridge University Press, Cambridge, 2003.

\bibitem[PS]{pisshlyak}
G.~Pisier, D.~Shlyakhtenko.
\newblock Grothendieck's theorem for operator spaces.
\newblock {\em Invent. Math.}, 150:185--217, 2002.

\bibitem[PX1]{px-BG}
G.~Pisier, Q.~Xu.
\newblock Non-commutative martingale inequalities.
\newblock {\em Comm. Math. Phys.}, 189:667--698, 1997.

\bibitem[PX2]{px-survey}
G.~Pisier, Q.~Xu.
\newblock Non-commutative {$L\sp p$}-spaces.
\newblock In {\em Handbook of the geometry of Banach spaces, Vol.\ 2}, pages
  1459--1517. North-Holland, Amsterdam, 2003.

\bibitem[R1]{ran-KP}
N.~Randrianantoanina.
\newblock Kadec-{P}e\l czy\'nski decomposition for {H}aagerup {$L\sp
  p$}-spaces.
\newblock {\em Math. Proc. Cambridge Philos. Soc.}, 132:137--154, 2002.

\bibitem[R2]{ran-mtrans}
N.~Randrianantoanina.
\newblock Non-commutative martingale transforms.
\newblock {\em J. Funct. Anal.}, 194:181--212, 2002.

\bibitem[R3]{ran-weak}
N.~Randrianantoanina.
\newblock A weak type inequality for non-commutative martingales and
  applications.
\newblock {\em Proc. London Math. Soc.}, 91:509--542, 2005.

\bibitem[R4]{ran-burk}
N.~Randrianantoanina.
\newblock Conditionned square functions for noncommutative martingales.
\newblock {\em Ann. Proba.}, to appear.

\bibitem[RX]{rx-JFA}
Y.~Raynaud, Q.~Xu.
\newblock On subspaces of non-commutative {$L\sb p$}-spaces.
\newblock {\em J. Funct. Anal.}, 203:149--196, 2003.

\bibitem[Ro]{ros-ineq}
H.~P.~Rosenthal.
\newblock On the subspaces of {$L\sp{p}$} {$(p>2)$} spanned by sequences of
  independent random variables.
\newblock {\em Israel J. Math.}, 8:273--303, 1970.


\bibitem[S]{shlya-quasifree}
D.~Shlyakhtenko.
\newblock Free quasi-free states.
\newblock {\em Pacific J. Math.}, 177:329--368, 1997.

\bibitem[St]{stratila}
{\c{S}}.~Str{\u{a}}til{\u{a}}.
\newblock {\em Modular theory in operator algebras}.
\newblock Editura Academiei Republicii Socialiste Rom\^ania, Bucharest, 1981.


\bibitem[Su]{suk-noniso}
F.~A.~Sukochev.
\newblock Non-isomorphism of {$L\sb p$}-spaces associated with finite and
  infinite von {N}eumann algebras.
\newblock {\em Proc. Amer. Math. Soc.}, 124:1517--1527, 1996.

\bibitem[T1]{tak-cond}
M.~Takesaki.
\newblock Conditional expectations in von {N}eumann algebras.
\newblock {\em J. Funct. Anal.}, 9:306--321, 1972.

\bibitem[T2]{tak-I}
M.~Takesaki.
\newblock {\em Theory of operator algebras. {I}}.
\newblock Springer-Verlag, New York, 1979.

\bibitem[VDN]{VDN}
D.~V.~Voiculescu, K.~J.~Dykema, A.~Nica.
\newblock {\em Free random variables}, volume~1 of {\em CRM Monograph Series}.
\newblock American Mathematical Society, Providence, RI, 1992.

\bibitem[X1]{xu-descrip}
Q.~Xu.
\newblock A description of {$(C\sb p[L\sb p(M)],R\sb p[L\sb p(M)])\sb \theta$}.
\newblock {\em Proc. Roy. Soc. Edinburgh Sect.}, 135:1073--1083, 2005.

\bibitem[X2]{xu-gro}
Q.~Xu.
\newblock Operator-space {G}rothendieck inequalities for noncommutative {$L\sb
  p$}-spaces.
\newblock {\em Duke Math. J.}, 131:525--574, 2006.

\bibitem[X3]{xu-embed}
Q.~Xu.
\newblock Embedding of {$C\sb q$} and {$R\sb q$} into noncommutative {$L\sb
  p$}-spaces, {$1\leq p<q\leq2$}.
\newblock {\em Math. Ann.}, 335:109--131, 2006.

\end{thebibliography}

\end{document}